\documentclass[11pt,a4paper]{article}
\usepackage[vmargin=24mm,hmargin=25mm]{geometry}
\usepackage{amsmath,amssymb,amsthm}
\allowdisplaybreaks[2]
\numberwithin{equation}{section}
\DeclareMathOperator{\Div}{div}

\newcommand{\R}{\mathbb R}
\newcommand{\dd}{\,\mathrm{d}x}
\newcommand{\dt}{\,\mathrm{d}t}
\newcommand{\eps}{\varepsilon}
\newtheorem{lemma}{Lemma}[section]
\newtheorem{theorem}{Theorem}[section]

\newtheorem{definition}{Definition}[section]

\newtheorem{proposition}{Proposition}[section]
\newcommand{\bnn}{\begin{eqnarray*}}
\newcommand{\enn}{\end{eqnarray*}}
\newcommand{\ba}{\begin{aligned}}
\newcommand{\ea}{\end{aligned}}
\newcommand{\be}{\begin{equation}}
\newcommand{\ee}{\end{equation}}
 
\newcommand{\la}{\label}

\begin{document}
 \medskip

\centerline{\Large\bf  Global large-data strong solutions to the 2D barotropic  }

\smallskip
\centerline{\Large\bf
 compressible Navier--Stokes equations with }

 \smallskip
\centerline{\Large\bf
linear bulk viscosity}

\bigskip

\centerline{Qian Li$^{1}$,\ \ Zhilei Liang$^{2}$,\ \ Zhuan Ye$^{3}$}

\medskip
\centerline{$^1$  School of Mathematics, China University of Mining and Technology,}
 \centerline{
Xuzhou 221116, Jiangsu, China}

\smallskip
\smallskip

\centerline{$^2$ School of Mathematics,
}

\centerline{Southwestern University of Finance and Economics, Chengdu 611130, P. R. China}

\smallskip
\smallskip

\centerline{$^3$ Department of Mathematics and Statistics, Jiangsu Normal University,
}

\centerline{101 Shanghai Road, Xuzhou 221116, Jiangsu, PR China}

\footnotetext{li-qian@mail.bnu.edu.cn (Qian Li), liangzl@swufe.edu.cn (Zhilei Liang), yezhuan815@126.com (Zhuan Ye)}

\bigskip

{\bf Abstract:}~~%
This paper concerns the global well-posedness of the two-dimensional (2D) compressible Navier--Stokes equations with density-dependent bulk viscosity and far-field vacuum. We consider the constitutive laws \(P(\rho)=\rho^\gamma\), with \(\gamma>1\), and \(\lambda(\rho)=\rho\), for which the bulk viscosity degenerates in vacuum.
We establish the existence and uniqueness of global strong solutions without any smallness restriction on the size of the initial energy. The central difficulty is to derive an upper density bound on arbitrary finite time intervals in the absence of enhanced bulk dissipation. The approach is to develop a logarithmic energy method based on interval layer foliation, which exploits the dissipation accumulated along particle trajectories as the density increases. Our result extends the large-data global strong-solution theory for density-dependent viscosity to the critical linear bulk-viscosity law in the whole two space dimensions.

{\vskip 1mm
 {\bf AMS Subject Classification 2020:}\quad 35Q35; 35B65; 76N10.

 }

\vskip .3in

\section{Introduction and main results}

We consider the Cauchy problem for the two-dimensional (2D)
barotropic compressible Navier--Stokes equations
\begin{equation}\label{eq:system}
\begin{cases}
 \rho_t+\Div(\rho u)=0,\\[2mm]
 (\rho u)_t+\Div(\rho u\otimes u)+\nabla P(\rho)
   =\mu\Delta u+\nabla\bigl((\mu+\lambda(\rho))\Div u\bigr),
\end{cases}
\end{equation}
where $\rho=\rho(x,t)\geq0$ and $u=(u_{1},u_{2})(x,t)$
denote the density and velocity, respectively. The pressure function is given by $P(\rho)=\rho^\gamma$ with the adiabatic exponent $\gamma>1$. Here,  the  shear viscosity $\mu$  is a constant  and   the bulk
viscosity $\lambda(\rho)$ is a linear function of density, in  particular,
\[
\mu>0,\qquad \lambda(\rho)=\rho.
\]
We supplement system \eqref{eq:system} with the far-field vacuum condition
in the appropriate sense, and prescribe the initial data
\begin{equation}\label{n4}
 \rho(x,t=0)=\rho_0(x),\qquad
  \rho u(x,t=0)=\rho_0(x)u_0(x),
 \quad x\in \mathbb R^2,
\end{equation}
where we define $u_0=0$ on the sets $\{ x\in \mathbb R^2: \rho_{0}(x)=0\}$.


\medskip

System \eqref{eq:system} describes the motion of a viscous
compressible fluid under a barotropic pressure law. The first
equation is the conservation of mass, and the second is the
balance of momentum. Density-dependent bulk viscosities are
important in the analysis of compressible flows because the
strength of the dissipation associated with $\Div u$ varies with
the density. This dependence creates additional difficulties in
vacuum regions, where the bulk-viscous contribution may degenerate.
Power-law constitutive relations of the form
$\lambda(\rho)=\rho^\beta$ provide a natural framework for
studying this phenomenon; see, for example,
\cite{Yang2006vacuum,YangYaoZhu2001,BreschDesjardinsLin2003,
MelletVasseur2007}.

For compressible Navier--Stokes equations with constant
viscosities, the global theory has been developed for several
classes of initial data and solutions. Global strong solutions
near equilibrium states were obtained by Matsumura and Nishida
\cite{MatsumuraNishida1980}. For large data, the foundational
works of Lions \cite{Lions1998} and Feireisl, Novotny, and
Petzeltova \cite{FeireislNovotnyPetzeltova2001} established the
existence of global finite-energy weak solutions. Recently, an important breakthrough was achieved by Hong, Peng, and Zhu \cite{zhu}, who extended  Lions's weak theory to $\gamma>4/3$, thereby covering the important case $\gamma=7/5$ for air. Results allowing
vacuum include the local classical theory of Cho and Kim
\cite{ChoKim2006} and the global classical solutions under a
suitable small-energy condition obtained by Huang, Li, and Xin
\cite{HuangLiXin2012}. These results concern different solution
classes and assumptions; in particular, they do not directly
settle the large-data strong-solution problem for the
density-dependent viscosity law considered here. We also mention
the finite-time implosion and singularity-formation results for
certain radially symmetric solutions of the three-dimensional
compressible Navier--Stokes and Euler equations due to Merle,
Raphael, Rodnianski, and Szeftel
\cite{MPRSAM1,MPRSAM2}.

The density-dependent bulk-viscosity model has a separate
development. For the constitutive law
$\lambda(\rho)=\rho^\beta$, Vaigant and Kazhikhov
\cite{Ka} established global strong solvability in two dimensions
for $\beta>3$ under assumptions that include a positive lower
bound for the density. Li and Liang \cite{hlma} developed the
local strong and classical solvability theory in the presence of
far-field vacuum. On the torus, Jiu, Wang, and Xin \cite{jwx}
treated initial densities that may vanish, and Huang and Li
\cite{hlia} subsequently obtained a global result in the range
$\beta>4/3$. For the whole-space Cauchy problem with far-field
vacuum, Huang and Li \cite{huang-li-cauchy-published}
established global existence and uniqueness of strong and
classical solutions for $\beta>4/3$ and $\gamma>1$.
More recently, Jiu and Zhang \cite{JiuZhang2026} obtained global
large radially symmetric strong solutions for $\beta\geq1$.
The radial symmetry in that result is an essential part of the
argument and does not cover general, non-radially symmetric data.

There are also related results under different hypotheses or on
different spatial domains. Jiu, Wang, and Xin studied the Cauchy
problem in weighted spaces \cite{jwx-weighted}. Zhang and Fang
proved global existence of weak solutions under a small-energy
assumption \cite{zf}, while Fang and Guo obtained global strong
solutions in the range $0\leq\beta\leq1$ under a smallness
condition on the total initial energy \cite{fg}. Chen, Zhang, and Zhu \cite{CZ1,CZ2} established the global
existence and uniqueness of classical solutions for large spherically symmetric initial
data with either vacuum or non-vacuum far-field density.  In a bounded
domain, Ducomet and Ne\v{c}asov\'a considered the
compressible Navier--Stokes--Fourier system with
density-dependent viscosities \cite{dn}. Huang, Su, Yan, and Yu
proved global existence of large radially symmetric strong
solutions in two-dimensional balls for a range of exponents
with $\beta>1$ \cite{hsyy}. These works provide important
results, but their assumptions on the size or symmetry of the
data, or their spatial setting, differ from the general
whole-space large-data problem addressed in this paper.

A central issue in the global theory is to obtain an
upper bound for the density that is uniform on every finite
time interval. Such a bound is particularly delicate in the
presence of far-field vacuum: the bulk-viscous dissipation
degenerates where $\rho$ is small, and the usual estimates must
also control the interaction among the effective viscous flux,
the vorticity, and the material acceleration. Recently, Fan, Jiang, Li, and L\"u \cite{FJLL} studied the periodic problem and the whole-space Cauchy
problem for the two-dimensional compressible Navier--Stokes equations with general density-dependent bulk viscosity. They established the global existence of strong solutions for large initial data allowing
vacuum in both settings. In particular, their results cover the superlinear power-law case
\[
 \lambda(\rho)=\rho^\beta,\qquad \beta>1,
\]
as well as the logarithmically enhanced linear law
\[
 \lambda(\rho)=\rho[\log(e+\rho)]^\alpha,
 \qquad \alpha>\frac12.
\]
The purpose of the present paper is to establish the corresponding
large-data global strong-solution result for the unenhanced
linear law $\lambda(\rho)=\rho$, without imposing radial
symmetry or a smallness condition on the total initial energy.
The key point is a logarithmic energy estimate at successive
density levels, combined with estimates for the effective
viscous flux and the material derivative. These estimates yield
the required finite-time density bound and allow the local
strong solution to be continued globally.

\bigskip

The aim of this paper is to establish the global existence
of strong solutions to problem \eqref{eq:system}--\eqref{n4}, without any smallness assumption on the initial data.
By strong solution, we mean
\begin{definition}[Strong solution]
For given $T>0$, the function  pair $(\rho,u)$ with $\rho\ge0$ is called a \emph{strong solution}
of \eqref{eq:system}--\eqref{n4} on $[0,T)$ if
  the equations in
\eqref{eq:system}  are satisfied a.e. in $\mathbb R^2\times (0,T)$, and
  the initial data are attained in the sense of distributions.
\end{definition}
Throughout the paper, we use the standard notation
\[ \int f(x)\,\dd=\int_{\mathbb{R}^2} f(x)\,\dd,\qquad \nabla^{\perp}=(\partial_{2}, -\partial_{1}),
\]
\[
L^p=L^p(\mathbb R^2),\qquad
W^{s,p}=W^{s,p}(\mathbb R^2),\qquad
H^s=W^{s,2}(\mathbb R^2),\qquad \|f\|_{p}=\|f\|_{L^p(\mathbb R^2)}.
\]
Our main result is stated in the following theorem.

\begin{theorem}[Global strong solution]\label{TH11}
Let $\gamma>1$, $q>2$, $a\in(1,2)$, and let
\begin{equation}\label{2.07}
 \bar x=(e+|x|^2)^{1/2}
          \log^{1+\eta_0}(e+|x|^2),\quad \eta_0>0.
\end{equation}
Suppose that the initial data satisfy
\begin{equation}
 \rho_0\geq0,\quad 0<M\triangleq \int \rho_0\,\dd,\quad \bar x^a\rho_0\in L^1\cap H^1\cap W^{1,q},\quad
 \nabla u_0\in L^2,\quad
 \sqrt{\rho_0}\,u_0\in L^2.\nonumber
\end{equation}
Then, for any $T>0,$ the Cauchy problem \eqref{eq:system}--\eqref{n4}
admits a unique global strong solution, with the following  properties
\begin{equation*}
\begin{aligned}
 &\rho\in C([0,T];L^1\cap H^1\cap W^{1,q}),\qquad
 \bar x^a\rho\in
 L^\infty(0,T;L^1\cap H^1\cap W^{1,q}),\\
 &\sqrt\rho\,u,\ \nabla u,\ \bar x^{-1}u,\
 \sqrt{t\rho}\,u_t
       \in L^\infty(0,T;L^2),\\
 &\nabla u\in L^2(0,T;H^1)
       \cap L^{(q+1)/q}(0,T;W^{1,q}),\\
 &\sqrt t\,\nabla u\in L^2(0,T;W^{1,q}),\qquad
 \sqrt\rho\,u_t,\ \sqrt t\,\nabla u_t,\
 \sqrt t\,\bar x^{-1}u_t
       \in L^2(0,T;L^2).
\end{aligned}
\end{equation*}
Moreover, there exists a constant $N>0$, depending only on $T$ and the initial data, such that
\begin{equation}
 \inf_{0\leq t\leq T}
 \int_{B_N}\rho(t,x)\,\dd
 \geq \frac14 \int \rho_0(x)\,\dd,\nonumber
\end{equation}
where $
 B_N=\{x\in\mathbb R^2:|x|<N\}
$ is a ball in $\mathbb R^2$ centered at the origin with radius $N$.
\end{theorem}

\bigskip
\textbf{Methodology.}
A key ingredient in the proof
is to obtain an upper bound for the density. However, the critical linear-viscosity case $\lambda(\rho)=\rho$ does not provide the
logarithmically enhanced dissipation. Our approach uses
 a density-level argument that combines weighted
energy estimates, a logarithmic estimate for the material derivative,
and a Lagrangian argument for the effective viscous flux.
Let us  briefly describe the contributions and new ideas in the proof in the following paragraphs.

\medskip

For any given $S<\infty,$  we  define
a finite but $S$ dependent number
\begin{equation}
 R=\max\left\{2,\sup_{0\leq t\leq S}\|\rho(t)\|_{L^\infty}\right\}.\nonumber
\end{equation}
The basic energy identity gives
$
 \int_0^t D(s)\,\mathrm ds=E(0)-E(t),
$
where the total energy $E$ and the   dissipation rate $D$ are given by
\[
 E(t)= \int
       \left(\frac12\rho|u|^2+\frac{P}{\gamma-1}\right)\dd,\quad D(t)= \int
 \bigl(\mu|\omega|^2+h(\rho)|\Div u|^2\bigr)\dd,
 \ \  h(\rho)=2\mu+\rho.
\]
By  spatial moment estimate and   the weighted Poincar\'e
inequality,  we deduce that
\[
 \|\sqrt{\rho}\,u\|_{{2}}^2
 \leq C\log(e+t)D(t).
\]
This estimate is  independent of $R$.

We  derive decay estimates for the kinetic energy and weighted estimates for the dissipation. We first derive a pressure
identity by testing the momentum equation with a potential generated
by a general function of $\rho$.  A bounded choice of the test
function gives the weighted dissipation estimate
\[
 \int_0^S\sqrt{\log(e+t)}\,D(t)\,\mathrm dt\leq C.
\]
To estimate the density moments, we introduce the transported
potential
\[
 \theta(\rho)=\rho-1+2\mu\log\rho,
 \qquad
 D_t\theta(\rho)+P=-G,
\]
where $G=h(\rho)\Div u-P$
is the effective viscous flux. Using this identity, together with a Riesz-transform commutator estimate and a weighted Trudinger inequality,
we are able to prove that, for all $p<\infty$,
\[
 \sup_{0\leq t\leq S}\|\rho(t)\|_{{p}}^p
 \leq C_{p,\varepsilon}R^\varepsilon,\quad
\]
Moreover,  by choosing the test function $g(\rho)=P(\rho)$ in the  pressure identity (see Lemma \ref{lem:pressure} below),  we  also have  \[
 \sup_{t\leq S}\|P(t)\|_{{2}}^2
 +\int_0^S\|P(t)\|_{{2}}^2\,\mathrm dt
 \leq C_\varepsilon R^\varepsilon.
\]
The important feature is that  we can choose  $\varepsilon>0$ arbitrarily small to control the  dependence on the  polynomial growth $R^{\varepsilon}$.

Next,  we organize the possible growth of the density into dyadic
levels.  Choose $r_j=2^jr_0$, and let $t_j$ be the first time at
which $\|\rho(t_j)\|_\infty=r_j$.  On the interval
$I_j=[t_{j-1},t_j]$, we define
\[
 E_j=\int_{I_j}D(t)\,\mathrm dt,\qquad
 Q_j=r_j^{-1/3}+\sum_{i=1}^j2^{i-j}E_i,\qquad
 L_j=\log(e+r_j).
\]
The intervals $I_j$ have disjoint interiors.  Consequently, the
basic energy inequality guarantees
\[
 \sum_j E_j\leq C,\qquad \sum_j Q_j\leq C.
\]
These  summation  quantities measure the dissipation spent in reaching
the successive density levels.

The central analytic estimate is a logarithmic energy inequality
for the material derivative.  For each fixed terminal level $r_j$,
we use the normalization
\[
 H_j(t)^2=e+r_j^\eta+A(t)^2,
 \qquad
 A(t)^2= \int
 \left(\mu|\omega|^2+\frac{|G|^2}{h(\rho)}\right)\,\mathrm dx,
\]
where $0<\eta<\min\{1,\gamma-1\}$.  The normalization is fixed throughout
the whole interval $[\tau,t_j]$; in particular, it is not allowed
to vary with time.  Multiplying  the momentum equation by  $\dot u$,  using  the  interpolation  and the pressure bounds, we obtain
\[
 \int_{I_j}\frac{B(t)^2}{H_j(t)^2}\,\mathrm dt
 \leq C_\tau r_jL_jQ_j,
 \qquad
 B(t)=\|\sqrt{\rho}\,\dot u\|_{{2}}.
\]
The fixed normalization is important: it allows the  integration
of the logarithmic differential inequality without differentiating a
time-dependent density-level parameter.

We then convert the preceding estimate into an $L^\infty$ bound for
the effective viscous flux on each density-crossing interval.  The
material-derivative estimate \eqref{sdbnmygf}
and  the elliptic relation
$
 \operatorname{div}(\rho\dot u)=\Delta G
$
and a potential interpolation inequality are used.  Choosing the interpolation
exponent with respect to  $L_j=\log(e+r_j)$ helps to  control  the   power growth
generated by the interpolation uniformly in $j$.  The result is
\[
 \int_{I_j}\frac{\|G(t)\|_\infty^2}{H_j(t)^2}\,\mathrm dt
 \leq C_\tau r_j^2L_j^2Q_j.
\]

The final step is a Lagrangian argument.  We consider  a particle path  along
which the density rises from $r_{j-1}=r_j/2$ to $r_j$.  Since $P=\rho^\gamma\geq 2^{-\gamma}r_j^\gamma$ on the part of the
trajectory, and  the positivity
produces a favorable  contribution in the identity
$
 D_t\theta(\rho)=-G-P,
$
  we use the   increase
of $\theta(\rho)$  and  the effective viscous flux estimate and the definition
of $H_j$, to conclude
\[
 \frac{1}{L_j}\leq C_\tau\sqrt{E_jQ_j}.
\]
This is the cost of crossing the $j$-th density level.  Summing over
all attained levels and using Cauchy--Schwarz yields
\[
 \sum_{j=1}^N\frac{1}{L_j}
 \leq C_\tau
 \left(\sum_{j=1}^N E_j\right)^{1/2}
 \left(\sum_{j=1}^N Q_j\right)^{1/2}
 \leq C_\tau.
\]
On the other hand, $r_j=2^jr_0$ implies $L_j\asymp j+1$, so that
\[
 \sum_{j=1}^\infty\frac{1}{L_j}
 \asymp\sum_{j=1}^\infty\frac{1}{j+1}=\infty.
\]
Thus, infinitely  many density levels are impossible.  The
contradiction guarantees that  the density is bounded from above,  independent of the terminal time $S$.

Once the density is bounded, the rest  regularity  estimates are standard. 
Roughly,  a  time-weighted $L^{2}$ norm estimate for  $\sqrt{t}\,\nabla\dot u$ as well as  the other high order derivatives
quantities are available.  
Having  these  uniform bounds of solutions  in hand,  we are allowed  to extend the local solution obtained in Lemma \ref{lem2.1} by restart new beginning time.
This rules out finite-time breakdown and yields the unique global existence of  strong
solution.  \\

The remainder of the paper is organized as follows. In Section~\ref{se:preliminaries}, we introduce the notations and provide some   basic   conservation laws, useful lemmas. In Section~\ref{Secapriori}, we derive the a priori estimates needed to control the solution before a uniform upper bound for the density is available. In particular, we establish the weighted energy and dissipation estimates, derive the pressure estimates, and develop the logarithmic energy method at successive density levels. Section~4 is devoted to proving a uniform upper bound for the density by combining the effective viscous flux and vorticity estimates with a Lagrangian argument along particle trajectories. In Section~\ref{sec:continuation}, we complete the proof of the main theorem. Finally, in Appendix \ref{sec:energy-tests}, we derive the weighted estimates for the material derivative.

\section{Preliminaries}\label{se:preliminaries}

In this section, we introduce the notation and establish several basic estimates that will be used throughout the paper.
We first provide  the local in time existence result of strong solution, whose proof
can be found in \cite{hlma}.
\begin{lemma}[Local existence]\la{lem2.1}
Under the same hypotheses listed in Theorem \ref{TH11},  there is a (small) time point $T_{0}>0$ which
 depends on $\mu$,$\gamma$,
$q,a,\eta_0$, $M$, and the initial values, such that the problem \eqref{eq:system}--\eqref{n4} has a unique strong solution over $(0,T_{0}]$.
\end{lemma}

\bigskip
For simplicity, the curl of the velocity field, the bulk viscosity and the effective viscous flux are separately denoted by
\[
 \omega=\nabla^{\perp}\cdot u\equiv\partial_2u_1-\partial_1u_2,\qquad
 h(\rho)=2\mu+\rho,\qquad  G=h(\rho)\Div u-P(\rho).
\]
In addition, the following  abbreviations are used:
\[
 \begin{aligned}
 A^2(t)&= \int
       \left(\mu|\omega|^2+\frac{|G|^2}{h(\rho)}\right)\dd,\\
 D(t)&= \int
       \left(\mu|\omega|^2+h(\rho)|\Div u|^2\right)\dd,\\
   E(t)&= \int
       \left(\frac12\rho|u|^2+\frac{P}{\gamma-1}\right)\dd,\\
 B(t)&=\|\sqrt\rho\,\dot u\|_{{2}},\qquad
 W(t)=\|\nabla\dot u\|_{{2}},\qquad \dot u=D_tu=\partial_{t}u+u\cdot\nabla u.
 \end{aligned}
\]
For every time \(t>0\) for which the solution exists,  $(\rho,u)$  of  \eqref{eq:system}--\eqref{n4} satisfies the conservation laws:
\be\la{ly2.2}
  \int \rho(t,x)\dd=M,\quad  \int \rho u(t,x)\dd
 = \int \rho_0u_0 (x)\dd,\quad
 E(t)+\int_0^t D(s)\mathrm{d} s =E(0).
\ee
We shall also use the function
\be\la{ly2.1}
 \theta(r)=\int_1^r\frac{h(z)}z\,dz\equiv r-1+2\mu\log r,\qquad r>0.
\ee
which implies
$
 \theta'(r)=\frac{h(r)}r>0,
$
and clearly, 
\begin{equation}\label{thaq211}
 D_t\theta(\rho)+P=-G.
\end{equation}

\bigskip

The following lemma concerns a two-dimensional  Poincar\'e type
inequality.
\begin{lemma}\la{lem2.2}(\cite[Lemma 3.3]{FJLL})
Let $\sqrt{\rho} u\in L^2(\mathbb R^2)$, $\nabla u\in L^2(\mathbb R^2)$ and
\[
 0\le \rho\in L^1(\mathbb R^2)\cap L^\gamma(\mathbb R^2),
 \qquad  \int \rho\dd=M>0,\quad
  \int \rho\log(e+|x|)\dd<\infty.
\]
Then,  the inequality  holds true
\[
  \int \rho\left|u-\frac1M \int \rho u\,\dd\right|^2\,\dd
 \leq C_\gamma
 \left(
 \|\rho\|_{\gamma}
 + \int \rho\log(e+|x|)\,\dd
 \right)\|\nabla u\|_{{2}}^2.
\]
\end{lemma}

\bigskip

In this paper, we may work in the zero-momentum Galilean frame.  Indeed, conservation of momentum implies that, for all $t>0$,
\[
M=\int\rho(x,t) \dd,\quad  \mathfrak m \triangleq\int \rho u(x,t) \dd
\]
are constant.
If we replace \(u\) with \(u-\mathfrak m/M\) and translate the spatial coordinates accordingly, then system \eqref{eq:system} remains invariant. Consequently, we may assume that
\be\la{ly6}
  \int \rho u(t,x)\dd=0.
\ee

We next deduce some space moment estimate on the kinetic energy.
By \eqref{2.07}, we see that
\[
 \left|\nabla\bigl(\bar x^{a/2}\bigr)\right|\leq C_{a,\eta_0}
 \qquad a\in (1,2).
\]
From the mass equation $\eqref{eq:system}_{1}$ we deduce that
\[
 \begin{aligned}
\frac{d}{dt} \int \bar x^a \rho(t,x) \dd
 &=2\int\rho\,\bar x^{a/2}u\cdot
                   \nabla(\bar x^{a/2})\,\dd\\
 &\leq C
       \left(\int\rho\bar x^a\,\dd\right)^{1/2}
       \left(\int\rho|u|^2\,\dd\right)^{1/2}\\
 &\leq C\left( \int \bar x^a \rho(t,x) \dd\right)^{1/2}E(t)^{1/2}
 \leq C \left( \int \bar x^a \rho(t,x) \dd\right)^{1/2},
 \end{aligned}
\]
where the last inequality follows from the bound of $E(t)\le E(0)$.
Thus, we obtain
\begin{equation}\label{eq:spatial-moment}
  \int \rho(t,x)\bar x^a\,\dd
 \leq C(e+t)^2.
\end{equation}
Letting $\phi(s)=\log(1+s)$, $f(x)=\bar x^a$ and noticing $\rho\geq0$, it follows from Jensen's inequality that
\[\phi\left(\frac{\int \rho(x,t) f(x) \dd}{\int \rho(x,t)  \dd}\right)\geq \frac{\int \rho(x,t) \phi(f(x)) \dd}{\int \rho(x,t)  \dd},
\]
which along with $M=\int\rho(x,t) \dd$ yields
\begin{align}\label{asdvn09}
 \int
\rho(t,x)\log(1+\bar x^a)\,\mathrm dx \leq
M\log\left(
1+\frac1M\int
\rho(t,x)\bar x^a\,\mathrm dx
\right).
\end{align}
Keeping in mind the fact
\[
\log(e+|x|)
\leq C_a+\log(1+\bar x^a),
\]
it then follows from \eqref{eq:spatial-moment} and \eqref{asdvn09} that
\[
\begin{aligned}
\int
\rho(t,x)\log(e+|x|)\dd
 \leq C_aM+
M\log\left(1+\frac{C}{M}(e+t)^2\right) \leq C\log(e+t).
\end{aligned}
\]
By this, \eqref{ly6}, and the fact $\|\rho(t)\|_{\gamma}^\gamma
 \leq  CE(0)$,   from Lemma \ref{lem2.2}  we conclude that
\be\la{ly7}
 \|\sqrt{\rho}\,u\|_{{2}}^2
 \leq C\log(e+t)\|\nabla u\|_{{2}}^2.
\ee

\section{The a priori estimates}
\label{Secapriori}
The local existence result in Lemma  \ref{lem2.1} guarantees that the problem admits a  unique strong solution $(\rho,u)$ in some short time interval. In the following, we
assume that $(\rho,u)$ exists on $[0,T]$ and deduce
the uniform a priori estimates under the same assumptions of Theorem \ref{TH11},
so that the local strong solution can be extended.

\medskip
The first lemma provides estimates for the kinetic
energy and the weighted dissipation.

\begin{lemma}\label{lem:weighted-dissipation}
For $0\leq t\leq S$, we have
\begin{equation}\label{eq:weighted-dissipation}
 E(t)\leq C\left(\frac{\ell(t)}{1+t}\right)^{(\gamma-1)/\gamma},
 \qquad \int_0^S\sqrt{\ell(t)}\,D(t)\dt\leq C,
\end{equation}
where $\ell(t)=\log(e+t)$.
\end{lemma}

\begin{proof}
Let
\be\la{ly2.3}
 I=(-\Delta)^{-1},\qquad
 T_{ij}=\partial_i\partial_j I,\qquad
 \psi=I\Div(\rho u).
\ee
It follows from the Hardy--Littlewood--Sobolev inequality that
\[
 \|\psi\|_{{2\gamma}}
 \leq C\|\rho u\|_{{2\gamma/(\gamma+1)}}
 \leq \|\rho\|_{{\gamma}}^{1/2}\|\sqrt{\rho}\,u\|_{ 2 }.
\]
Taking  divergence on  momentum equation yields\[
 \partial_t\Div(\rho u)
 +\partial_i\partial_j(\rho u_i u_j)
 =\Delta G,
\]
Operating it by $I=(-\Delta)^{-1}$ and using \eqref{ly2.3}, we have
\begin{equation}\label{eq:psi-evolution}
 \psi_t+G+ T_{ij}(\rho u_i u_j)=0.
\end{equation}

Denote by
\[
 b_g(r)=rg'(r)-g(r),\quad
 \mathcal H_g(r)
   =r\int_0^r\frac{h(z)g(z)}{z^2}\,\mathrm dz,
 \quad
 J_g(t)= \int \rho u\cdot\nabla Ig\dd,
\] where the function $g(\cdot)$ will be specified later.
It follows from the  continuity equation that
\[
 \begin{aligned}
 g_t+\Div(gu)
 =g'(\rho)\rho_t+\Div(g(\rho)u) =\bigl(g(\rho)-\rho g'(\rho)\bigr)\Div u =-b_g(\rho)\Div u,
 \end{aligned}
\]
and hence,
\[
 -\int\psi g_t\,\mathrm dx
 =-\int\psi\bigl(-\Div(gu)-b_g \Div u\bigr)\,\mathrm dx
 =-\int gu\cdot\nabla\psi+\int\psi b_g \Div u\,\mathrm dx.
\]
Since $
 \mathcal H_g'(r)
 =\int_0^r\frac{h(z)g(z)}{z^2}\,\mathrm dz
       +\frac{h(r)g(r)}r,
$
we have
\[
 \begin{aligned}
 \frac{\mathrm d}{\mathrm dt}\int\mathcal H_g(\rho)\,\mathrm dx
  =\int\bigl(\mathcal H_g(\rho)
        -\rho\mathcal H_g'(\rho)\bigr)\Div u\,\mathrm dx =-\int h(\rho)g(\rho)\Div u\dd,
 \end{aligned}
\]
where we have used
$\rho\mathcal H_g'(\rho)-\mathcal H_g(\rho)=h(\rho)g(\rho)$.
Observe from the symmetry of $T_{ij}$   that
\[
 \begin{aligned}
  \int T_{ij}(\rho u_i u_j)g\,\mathrm dx
       -\int gu\cdot\nabla\psi\,\mathrm dx
 & =
       \int\rho u_i u_jT_{ij}g\,\mathrm dx
       -\int\rho u_iT_{ij}(u_jg) \,\mathrm dx = \int\rho u_i[u_j,T_{ij}]g\,\mathrm dx.
 \end{aligned}\]
Therefore, taking the above identities into account,
we multiply \eqref{eq:psi-evolution} by $g$ and obtain
\[
 \begin{aligned}
 \left(J_g+\int\mathcal H_g(\rho)\,\mathrm dx\right)'
 &=\int (G-h\Div u)g\,\mathrm dx + \int\rho u_i[u_j,T_{ij}]g\,\mathrm dx
       +\int\psi b_g(\rho)\Div u\dd,
 \end{aligned}
\] where $$
 J_g(t)= \int \rho u\cdot\nabla Ig\dd=- \int \psi g\,\mathrm dx,$$ which follows from  \eqref{ly2.3} and  integration by parts.
Equivalently,
\begin{equation}\label{eq:pressure-identity}
 \left(J_g+\int\mathcal H_g(\rho)\dd\right)'
       +\int Pg\,\mathrm dx
 = \int \rho u_i[u_j,T_{ij}]g\dd
       +\int\psi b_g(\rho)\Div u\dd.
\end{equation}

Let us  estimate the right-hand side terms one by one. Set
\[
 n=\frac{2\gamma}{\gamma-1},\qquad
 \widetilde{q}=\frac{2\gamma}{\gamma+1},
 \qquad  \frac{1}{\widetilde{q}}+\frac1n=1.
\]
The Coifman--Rochberg--Weiss commutator estimate (cf.\cite{coi1,coi3}) and the fact
$\|u\|_{\mathrm{BMO}}\leq C\|\nabla u\|_{ 2 }$ shows that
\[
 \|[u_j,T_{ij}]g\|_{ n }
 \leq C\|\nabla u\|_{ 2 }\|g\|_{ n }.
\]
Therefore, by H\"older's inequality,
\[
 \begin{aligned}
 \left| \int\rho u_i[u_j,T_{ij}]g\right|
 &\leq C\|\rho u\|_{\widetilde{q}}\|\nabla u\|_{ 2 }\|g\|_{ n }.
 \end{aligned}
\]
For the other term, since
$
 \frac1{2\gamma}+\frac1n+\frac12=1,
$
 one gets
\[
 \left|\int\psi b_g(\rho)\Div u\,\mathrm dx\right|
 \leq \|\psi\|_{2\gamma}\|b_g(\rho)\|_{ n }\|\Div u\|_{ 2 }.
\]
Thanks to \eqref{ly7} and inequality $\|\psi\|_{{2\gamma}}\leq C\|\rho u\|_{\widetilde{q}}$,
\[
 \begin{aligned}
 \|\rho u\|_{\widetilde{q}}+\|\psi\|_{{2\gamma}}
 \leq C\|\rho\|_{{\gamma}}^{1/2}\|\sqrt\rho\,u\|_{ 2 } \leq C\sqrt{\ell(t)D(t)}.
 \end{aligned}
\]
Keeping in mind
$\|\Div u\|_{ 2 }\leq C\|\nabla u\|_{ 2 }\leq C\sqrt D$,
we derive that
\begin{equation}\label{eq:pressure-source}
 \left|\text{right-hand side of }\eqref{eq:pressure-identity}\right|
 \leq C\bigl(\|g\|_{ n }+\|b_g(\rho)\|_{ n }\bigr)
                \sqrt{\ell(t)}D(t).
\end{equation}
Let us  choose   function $g$ as
\be
 g_0(r)=\frac{P(r)}{h(r)+P(r)}
       =\frac{r^\gamma}{2\mu+r+r^\gamma}.\nonumber
\ee
Then, we easily check
\be\la{ly4.1} 0\leq g_0\leq1,\quad
 0\leq g_0(r)\leq \frac{P(r)}{2\mu},
 \quad |b_{g_0}(r)|\leq C.
\ee
The last inequality in \eqref{ly4.1} can be explained as below:
since \[
 rg_0'(r)
 =g_0(r)(1-g_0(r))
       \left(\gamma-\frac{rh'(r)}{h(r)}\right),
\] and the   $rh'(r)/h(r)\equiv\frac{r}{2\mu+r}$ is bounded,
 we have
\[
 \frac{b_{g_0}(r)}{g_0(r)}
 =(1-g_0(r))\left(\gamma-\frac{rh'(r)}{h(r)}\right)-1
\]
is bounded.  By \eqref{ly4.1},  there holds
\[
 \|g_0\|_{ n }^n+\|b_{g_0}(\rho)\|_{ n }^n
 \leq C\int g_0\,\mathrm dx
 \leq C\|P\|_{{1}}\leq C.
\]

For the left hand side of \eqref{eq:pressure-identity},  we see that
 $\mathcal H_{g_0}$ is nonnegative and
satisfies
\[
 \begin{aligned}
 \mathcal H_{g_0}(r)
 &=r\int_0^r\frac{h(z)}{z^2}
              \frac{P(z)}{h(z)+P(z)}\,\mathrm dz \leq r\int_0^r z^{\gamma-2}\,\mathrm dz
 =\frac{r^{\gamma}}{\gamma-1}
 =\frac{P(r)}{\gamma-1}.
 \end{aligned}
\]
Next  to bound $J_{g_0}$. By the  interpolation inequality,
\begin{equation}\label{eq:LP-interpolation13}
\|\nabla I f\|_{ \infty }= \|\nabla (-\Delta)^{-1} f\|_{ \infty }
    \leq C
    \|f\|_{ 1 }^{1-\frac{p}{2(p-1)}}
    \|f\|_{ p }^{\frac{p}{2(p-1)}},\quad \forall\,p\in(2,\infty].
\end{equation}
Substituting $f$ with $g_{0}$ in  \eqref{eq:LP-interpolation13} and choosing  $p=\infty$, we find
\[\|\nabla Ig_0\|_{ \infty }=
 \|\nabla (-\Delta)^{-1} g_0\|_{ \infty }
 \leq C\|g_0\|_{ 1 }^{1/2}\|g_0\|_{ \infty }^{1/2}
 \leq C\|P\|_{ 1 }^{1/2}.
\]
Consequently,
\[
 \begin{aligned}
 |J_{g_0}|
 &\leq \|\sqrt\rho\,u\|_{ 2 }
       \|\sqrt\rho\|_{ 2 }\|\nabla Ig_0\|_{ \infty } \leq C E^{1/2}\|P\|_{ 1 }^{1/2}
 \leq CE.
 \end{aligned}
\]
Since $h(r)\leq C(1+r^\gamma)$, we compute
\[
 \begin{aligned}
 \int\rho(h+P)^{(\gamma-1)/\gamma}\,\mathrm dx
 &\leq C\int\rho(1+\rho^\gamma)^{(\gamma-1)/\gamma}
                  \dd \leq C\int(\rho+\rho^\gamma)\dd
 \leq C.
 \end{aligned}
\]
Let
$
 a_0=\frac{\gamma-1}{2\gamma-1},
b_0=\frac{\gamma}{2\gamma-1}.
$
Then
$$
 P
 =\left(\frac{P^2}{h+P}\right)^{a_0}
   \left(\rho(h+P)^{(\gamma-1)/\gamma}\right)^{b_0}.
$$
Thus, by  H\"older's inequality,
\begin{equation}\label{eq:pressure-coercivity}
 \|P\|_{ 1 }
 \leq
 \left(\int Pg_0\,\mathrm dx\right)^{a_0}
 \left(\int\rho(h+P)^{(\gamma-1)/\gamma}\,\mathrm dx
       \right)^{b_0}
 \leq C\left(\int Pg_0\,\mathrm dx\right)^{a_0}.
\end{equation}
We define
\[
 \xi=\frac1{a_0}=\frac{2\gamma-1}{\gamma-1},
 \qquad
 \delta(t)=\frac{\delta_0}{\sqrt{\ell(t)}},
 \qquad
 \mathcal F(t)=E(t)+\delta(t)
       \left(J_{g_0}(t)+\int\mathcal H_{g_0}(\rho(t))\,\mathrm dx\right).
\]
The bounds for $J_{g_0}$ and $\mathcal H_{g_0}$ show that
\[
 \left|J_{g_0}+\int\mathcal H_{g_0}(\rho)\,\mathrm dx\right|
 \leq CE.
\]
Choosing $\delta_0>0$ sufficiently small, we obtain
\begin{equation}\label{eq:F-equivalence}
 cE(t)\leq\mathcal F(t)\leq CE(t).
\end{equation}
Since $E'=-D$, we deduce from \eqref{eq:pressure-identity} and
\eqref{eq:pressure-source}  that
\[
 \mathcal F'+D+\delta\int Pg_0\,\mathrm dx
 \leq C\delta\sqrt\ell D
       +|\delta'|\left|J_{g_0}+\int\mathcal H_{g_0}\,\mathrm dx\right|.
\]
Because $\delta\sqrt\ell=\delta_0$, the term on the right-hand side $C\delta\sqrt\ell D$ can be absorbed  by  $D$ on the left  by decreasing $\delta_0$ once more.
Using \eqref{eq:F-equivalence}, we arrive at
\begin{equation}\label{eq:F-differential}
 \mathcal F'+cD+\delta\int Pg_0\,\mathrm dx
 \leq C|\delta'|\mathcal F.
\end{equation}

For the kinetic energy, since it is bounded and obeys
$ \|\sqrt\rho\,u\|_{ 2 }^2
       \leq C\ell(t)D(t),
$
we thus obtain for $\xi>1$,
\[
\left(\|\sqrt\rho\,u\|_{ 2 }^2\right)^{\,\xi}\leq C\|\sqrt\rho\,u\|_{ 2 }^2
       \leq C\ell(t)D(t).
\]
For the internal energy, from \eqref{eq:pressure-coercivity} we get
\[
 \left(\frac1{\gamma-1}\|P\|_1\right)^\xi
 \leq C\int Pg_0\,\mathrm dx.
\]
Thanks to the last two inequalities and
\eqref{eq:F-equivalence}, we see  that
\begin{equation}\label{eq:F-coercive}
 \mathcal F^\xi
 \leq C\left(\ell(t)D(t)+\int Pg_0\,\mathrm dx\right).
\end{equation}
Since $\delta(t)=\delta_0\ell(t)^{-1/2}\geq
c\ell(t)^{-1}$, we deduce from
\eqref{eq:F-differential} and \eqref{eq:F-coercive} that
\begin{equation}\label{eq:F-bernoulli}
 \mathcal F'+c\ell(t)^{-1}\mathcal F^\xi
 \leq C|\delta'(t)|\mathcal F.
\end{equation}
To solve \eqref{eq:F-bernoulli}, we denote
\[
 Z(t)=\Phi(t)\mathcal F(t)\quad {\rm with}\quad \Phi(t)=\exp\left(-C\int_0^t|\delta'(s)|\,\mathrm ds\right),
\]
where the constant $C$ is the same as  in \eqref{eq:F-bernoulli}.
Since
$
 \int_0^\infty|\delta'(s)|\,\mathrm ds<\infty,
$
we thus obtain
$$c\leq\Phi(t)\leq1.$$
By this and \eqref{eq:F-bernoulli},  we differentiate $Z$ and get
\[
 Z'(t)+c\ell(t)^{-1}Z(t)^\xi\leq0,
\]
which can be explicitly solved
\[
 Z(t)\leq C\left(1+\int_0^t\frac{\mathrm ds}{\ell(s)}
            \right)^{-1/(\xi-1)}.
\]
Noting that
$
 \int_0^t\frac{\mathrm ds}{\log(e+s)}
 \geq c\,\frac{1+t}{\log(e+t)}
$ for all $t\ge1$, by a direct comparison, the same estimate also holds for sufficiently small positive times. For \(0\leq t\leq1\), the desired estimate follows
directly from the basic energy bound \(E(t)\leq E(0)\),
after enlarging the constant \(C\).
Therefore, we derive
\begin{equation*}
 E(t)\leq C\left(\frac{\ell(t)}{1+t}\right)^{(\gamma-1)/\gamma},\quad t\ge 0.
\end{equation*}

We finally derive the weighted dissipation bound.  Since $E'=-D$,
integration by parts gives
\[
 \begin{aligned}
 \int_0^S\sqrt{\ell(t)}D(t)\,\mathrm dt
 &=-\int_0^S\sqrt{\ell(t)}E'(t)\,\mathrm dt\\
 &=E(0)-\sqrt{\ell(S)}E(S)
       +\int_0^S E(t)\frac{\mathrm dt}
          {2(e+t)\sqrt{\ell(t)}}\\
 &\leq E(0)+C\int_0^\infty
       \frac{\ell(t)^{(\gamma-1)/\gamma-1/2}}
            {(1+t)^{(\gamma-1)/\gamma}(e+t)}
       \,\mathrm dt\\
 &\leq C,
 \end{aligned}
\]
where the last integral converges due to $(\gamma-1)/\gamma>0$.
Therefore, we arrive at
\begin{equation}
\int_0^S\sqrt{\ell(t)}D(t)\,\mathrm dt\leq C,\nonumber
\end{equation}
with a constant independent of $R$ and $S$.
\end{proof}

\bigskip

We next establish uniform-in-time \(L^p\) bounds for both the density and the  pressure needed
before the density upper bound is known. To this end, we
recall the a priori quantity
\[
R\triangleq \max\left\{
2,\sup_{0\leq t\leq S}\|\rho(t)\|_\infty
\right\}.
\]
\begin{lemma}\label{lem:pressure}
For every fixed  $p\in [1,\infty)$, there holds
\begin{equation}\label{eq:pressure1}
 \sup_{0\leq t\leq S}\|\rho(t)\|_p^p
      \leq C_{p,\eps}R^\eps,
\end{equation}
\begin{equation}\label{eq:pressure2}
 \sup_{0\leq t\leq S}\bigl(\|P(t)\|_2^2+\|P(t)\|_4^{4/3}\bigr)
                  +\int_0^S\|P(t)\|_2^2\dt\leq C_\eps R^\eps,
\end{equation}
where the number $\eps>0$ can be arbitrarily small.
\end{lemma}

\begin{proof}
Let $\psi(t,x)$ be as in  \eqref{ly2.3} and let $f_{+}=\max\{f,0\}$. Putting
\[
 M_p(t)= \int \rho \bigl(\theta(\rho(t,x))-\psi(t,x)\bigr)_+^{p-1}\dd,\quad (p\ge4),
\]
where    the integrand is assumed to be zero on vacuum sets.  


Making use of the embedding inequality, interpolation inequality,  and the  energy bound \eqref{ly2.2}, one has
\bnn\ba
 \|\psi\|_{ p }
 \leq C\|\rho u\|_{{2p/(p+2)}} \leq C\|\rho\|_{{p/2}}^{1/2}\|\sqrt\rho\,u\|_{ 2 }
 \leq C\|\rho\|_{{p/2}}^{1/2}\le CM^{\frac{1}{2(p-1)} }\|\rho\|_{ p }^{ \frac{p-2}{2(p-1)}}.
\ea\enn  Thus,
\begin{equation}\label{eq:psi-moment}
 \int\rho|\psi|^{p-1}\dd
 \leq \|\rho\|_{ p }\|\psi\|_{ p }^{p-1}
 \leq C\|\rho\|_{ p }^{\,1+(p-2)/2}
 =C\|\rho\|_{ p }^{p/2}.
\end{equation}
Notice that, by the definition of \(\theta\), for
\(\rho\geq2\), we have
\[
\rho\leq C\theta(\rho)
\leq C\left[
(\theta(\rho)-\psi)_++|\psi|
\right].
\]
Consequently,
\[
\rho^p
\leq C\rho\left[
(\theta(\rho)-\psi)_+^{p-1}
+|\psi|^{p-1}
\right].
\]
Integrating over \(\mathbb R^2\), we obtain
\begin{align}
\int \rho^p\,dx
&\leq C+C M_p(t)
+C\int \rho|\psi|^{p-1}\,dx
\nonumber\\
&\leq C+C M_p(t)
+C\left(\int \rho^p\,dx\right)^{1/2},\nonumber
\end{align}
where we have used \eqref{eq:psi-moment}.
By Young's inequality, it follows that
\begin{equation}\label{eq:moment-comparison}
\int \rho^p\,dx
\leq C_p(1+M_p(t)).
\end{equation}
Next, it follows from \eqref{thaq211} and
\eqref{eq:psi-evolution} that
\bnn
 D_t(\theta(\rho)-\psi)_{+}+P
 =- [u_i,T_{ij}](\rho u_j).
\enn
Multiplying this equation by
$(p-1)\rho \bigl(\theta(\rho(t,x))-\psi(t,x)\bigr)_+^{p-2}$ gives
\be\la{ly2.5}\ba
  & \frac{d}{dt}M_p
+(p-1)\int\rho P \bigl(\theta(\rho(t,x))-\psi(t,x)\bigr)_+^{p-2}\dd\\
 &
 =-(p-1)\int
       \rho  \bigl(\theta(\rho(t,x))-\psi(t,x)\bigr)_+^{p-2}[u_i,T_{ij}](\rho u_j)\dd.
 \ea\ee
For the last integral in \eqref{ly2.5}, applying H\"older's inequality and the commutator, we conclude
\[
 \begin{aligned}
 &\int
       \rho  \bigl(\theta(\rho(t,x))-\psi(t,x)\bigr)_+^{p-2}[u_i,T_{ij}](\rho u_j)\dd \\
       &\leq
 \left(\int\rho  \bigl(\theta(\rho(t,x))-\psi(t,x)\bigr)_+^{p-1}\right)^{(p-2)/(p-1)}
 \left(\int\rho
       \left| [u_i,T_{ij}](\rho u_j)\right|^{p-1}
       \right)^{1/(p-1)}\\
 & \leq
 C M_p^{(p-2)/(p-1)}\|\rho\|_{ p }^{1/(p-1)}
       \left\| [u_i,T_{ij}](\rho u_j)\right\|_{ p }\\
 &\leq C  M_p^{(p-2)/(p-1)}\|\rho\|_{ p }^{1/(p-1)}\|\nabla u\|_2\|\rho u\|_{ p }.
 \end{aligned}\]
Inserting it back into \eqref{ly2.5} and utilizing \eqref{eq:moment-comparison}, we get
 \begin{align}\label{eq:moment-ode}
 \frac{d}{dt}M_p
& +(p-1)\int\rho P \bigl(\theta(\rho(t,x))-\psi(t,x)\bigr)_+^{p-2}\dd\leq C_p(1+M_p)^{(p-1)/p}
          \|\nabla u\|_{ 2 }\|\rho u\|_{ p }.
 \end{align}

We are now ready to  estimate $\|\rho u\|_p$.  For this purpose, we present
  the weighted Trudinger inequality (cf.  \cite[Inequality (3.10)]{FJLL})
\begin{equation}\label{eq:weighted-trudinger}
  \int (e+t)^{-1}\bar x^{-4}
       \exp\left(\frac{\eps_{1}|u|^2}{\|\nabla u\|_{ 2 }^2}\right)
       \,\mathrm dx\leq C,
\end{equation} where $\eps_{1}>0$ is a small constant.

Let  $y,z\geq0$, then it holds
\be\la{ly2.7}
 yz^{p/2}
 \leq C_{p,\eps_{1}}\left(e^{\eps_{1}z}
       +y[\log(e+y)]^{p/2}\right)
\ee
(Hint: the inequality follows by separating $z$ into  $z\leq C\log(e+y)$ and
$z>C\log(e+y)$.)

Choose in  \eqref{ly2.7} the functions
\[y=(e+t)\bar x^4\rho^p,\qquad z=\frac{|u|^2}{\|\nabla u\|_{ 2 }^2}.
\]
Multiplying by $(e+t)^{-1}\bar x^{-4}$ and integrating the resulting expression,  we use \eqref{eq:weighted-trudinger} to deduce
 \begin{align}\label{eq:orlicz}
 \frac{\|\rho u\|_{ p }^p}{\|\nabla u\|_{ 2 }^p}
 &\leq C_p
 \left[
  \int (e+t)^{-1}\bar x^{-4}
       \exp\left(\frac{\eps_{1}|u|^2}{\|\nabla u\|_{ 2 }^2}\right)
       \,\mathrm dx+\int\rho^p
  \left\{\log\bigl(e+(e+t)\bar x^4\rho^p\bigr)\right\}^{p/2}
       \,\mathrm dx\right]\nonumber\\
        &\leq C_p
 \left[1+\int\rho^p
  \left\{\log\bigl(e+(e+t)\bar x^4\rho^p\bigr)\right\}^{p/2}
       \,\mathrm dx\right]\nonumber\\
   &\leq C_p
 \left[1+\int\rho^p
  \left\{\bigl(\ell(t)+\log(e+\rho)+\log\bar x\bigr)\right\}^{p/2}
       \,\mathrm dx\right].
 \end{align}

It remains to handle the spatial logarithm $\log\bar x$ in \eqref{eq:orlicz}. Setting
$
 Y=1+\frac{\bar x}{(e+t)^{2/a}}.
$
Then, from  \eqref{eq:spatial-moment} we obtain
\be\la{ly3.1}
\log\bar x\leq C_a\ell(t)+\log Y,\qquad  \int\rho Y^a\,\mathrm dx\leq C.
\ee
We  claim that
\begin{equation}\label{eq:spatial-log}
 \int\rho^p(\log Y)^{p/2}\dd
 \leq C_p\left(
 1+\int\rho^p[\log(e+\rho)]^{p/2}\,\mathrm dx\right).
\end{equation}
To prove this, we express
\[\int\rho^p(\log Y)^{p/2}\dd=\int_{\Omega}\rho^p(\log Y)^{p/2}\dd+\int_{\mathbb R^2\setminus\Omega}\rho^p(\log Y)^{p/2}\dd \]
with
 \[\Omega=\left\{x: \log Y\leq\frac{2(p-1)}{a}\log(e+\rho)\right\}.
\]
On the one hand,
\[\int_{\Omega}
 \rho^p(\log Y)^{p/2} \dd
 \leq C_p \int_{\Omega}\rho^p[\log(e+\rho)]^{p/2} \dd.
\]
On the other hand,  for $Y\geq1$  we have $(\log Y)^{p/2}\leq C_p Y^{a/2}$.  This together with   \eqref{ly3.1} implies
\[\int_{\mathbb R^2\setminus\Omega}\rho^p(\log Y)^{p/2}\dd
 =\int_{\mathbb R^2\setminus\Omega}\rho\rho^{p-1}(\log Y)^{p/2}\dd
 \leq C_p\int_{\mathbb R^2\setminus\Omega}\rho Y^a \dd \le C
\]
The desired \eqref{eq:spatial-log} thus follows.

In terms of
\eqref{eq:orlicz}--\eqref{eq:spatial-log} and
\eqref{eq:moment-comparison}, we get
\[
 \begin{aligned}
 \|\rho u\|_{ p }^p
 &\leq C_p\|\nabla u\|_2^p
       \left[1+\{\ell(t)\log(e+R)\}^{p/2}
                         (1+M_p)\right],
 \end{aligned}
\]
that is,
\begin{equation}\label{eq:rho-u-moment}
 \|\rho u\|_{ p }
 \leq C_p \|\nabla u\|_{ 2 }\sqrt{\ell(t)\log(e+R)}\,
       (1+M_p)^{1/p}.
\end{equation}
Inserting \eqref{eq:rho-u-moment} into \eqref{eq:moment-ode} yields
\begin{equation}\label{eq:moment-growth}
\frac{\mathrm d}{\mathrm dt}\log(1+M_p)
 \leq C_p\sqrt{\ell(t)\log(e+R)}\,\|\nabla u\|_{ 2 }^2
 \leq C_p\sqrt{\ell(t)\log(e+R)}\,D(t).
\end{equation}
Using \eqref{eq:weighted-dissipation},  we integrate   \eqref{eq:moment-growth} to obtain
\[
 \sup_{0\leq t\leq S}\log(1+M_p(t))
 \leq C_p+C_p\sqrt{\log(e+R)},
\]
which together  with \eqref{eq:moment-comparison}  implies that, for any small $\varepsilon>0$,
\bnn
 \sup_{0\leq t\leq S}\|\rho(t)\|_{ p }^{p}
 \leq C_p\exp\!\left(C_p\sqrt{\log(e+R)}\right) \leq C_{p,\varepsilon}R^\varepsilon.
\enn
Regarding $p\in [1,4)$, it follows that
\[
\|\rho(t)\|_p^{p}
\leq
\|\rho(t)\|_1^{\frac{4-p}{3}}
\|\rho(t)\|_4^{p-\frac{4-p}{3}}\leq C_{\varepsilon}R^\varepsilon.
\]
This proves   \eqref{eq:pressure1}.

Next, let us  pay attention to  \eqref{eq:pressure2}.
Considering  \eqref{eq:pressure-identity} with $g=P=\rho^\gamma$, we have
\[
 b_P(\rho)=\rho\frac{\mathrm d}{\mathrm d\rho}\rho^\gamma-\rho^\gamma
          =(\gamma-1)P,
\quad
 \mathcal H_P(r)
 =r\int_0^r h(z)z^{\gamma-2}\,\mathrm dz\geq0.
\]
Also, the same deduction as  \eqref{eq:pressure-source}  shows
\[
 \left|\text{right-hand side of }\eqref{eq:pressure-identity}\right|
 \leq C_\gamma\|P\|_{ n }\sqrt{\ell(t)}D(t).
\]
Recalling \eqref{eq:pressure1}, one has
\[
 \sup_{0\leq t\leq S}\|P\|_{ n }
 =\sup_{0\leq t\leq S}\|\rho\|_{{\gamma n}}^{\gamma}
 \leq C_\varepsilon R^\varepsilon,
\]
With the last two inequalities, as well as Lemma \ref{lem:weighted-dissipation},  we obtain after integration
\begin{equation}\label{eq:pressure-source-integrated}
 \int_0^S
 \left|\text{right-hand side of }\eqref{eq:pressure-identity}\right|
 \,\mathrm dt
 \leq C_\varepsilon R^\varepsilon
       \int_0^S\sqrt{\ell(t)}D(t)\,\mathrm dt
 \leq C_\varepsilon R^\varepsilon.
\end{equation}

Next, we estimate $J_P$. By \eqref{eq:LP-interpolation13} with $p=4$, we have
\[
 \begin{aligned}
 \|\nabla IP\|_{ \infty }
 \leq   C(\|P\|_{ 1 }+\|P\|_{{4}}).
 \end{aligned}
\]
Thus,
\begin{align}
 |J_P(t)|
 &\leq\|\rho u\|_{ 1 }\|\nabla IP\|_{ \infty }\nonumber\\&
 \leq C\|\sqrt\rho\,u\|_{ 2 }\|\sqrt\rho\|_{ 2 }(\|P\|_{ 1 }+\|P\|_{{4}})
 \nonumber\\&\leq C(1+\|P\|_{{4}}).\nonumber
\end{align}
Keeping in mind \eqref{eq:pressure1}, we have
\be\la{ly3.3}
 \sup_{0\leq t\leq S}\|P\|_{ 4 }
 \leq C_\varepsilon R^\varepsilon.
\ee
Consequently,    $|J_P(t)|\leq C_\varepsilon R^\varepsilon$.
Integrating \eqref{eq:pressure-identity}
with $g=P$,  using \eqref{eq:pressure-source-integrated}  and  \eqref{ly3.3}, we infer
\be\la{ly3.4}
 \begin{aligned}
 \int_0^t\|P(s)\|_{ 2 }^2\,\mathrm ds
 &=J_P(0)-J_P(t)
   +\int\mathcal H_P(\rho(0))\,\mathrm dx
   -\int\mathcal H_P(\rho(t))\,\mathrm dx\\
 &\quad+\int_0^t
       \left(\text{right-hand side of }
                    \eqref{eq:pressure-identity}\right)\mathrm ds\\
                    &\le C_\varepsilon R^\varepsilon,
 \ea\ee
where we have used the fact $\mathcal H_P(t)\geq0$. Therefore, the combination of \eqref{ly3.3} with \eqref{ly3.4} concludes \eqref{eq:pressure2}.
\end{proof}

We next derive an energy estimate on each interval during which the
density crosses one prescribed level. Before introducing the density levels \(r_j\), we first fix a
time \(\tau\in(0,T_0)\) such that
\[
A(\tau)+B(\tau)<\infty.
\]
Such a choice is possible because the local-in-time estimates
established above imply that \(A(t)\) and \(B(t)\) are finite
for almost every \(t\in(0,T_0)\). We fix this \(\tau\)
throughout the density-level argument. In particular, all
estimates below are carried out on the intervals
\([\tau,t_j]\), with \(t_0=\tau\), and the constants denoted
by \(C_\tau\) may depend on this fixed choice of \(\tau\), but
not on  the density level \(r_j\) and the index
\(j\).
To proceed, we fix
\[
 0<\eta<\min\{1,\gamma-1\},\qquad
 0<\varepsilon\leq\min\{\eta,1/12\}.
\]
We choose sufficiently large $r_0$ such that
$$r_0\geq 2+\sup_{0\leq t\leq T_0}\|\rho(t)\|_{\infty}.$$
In particular, none of the levels \(r_j\) $(j\geq1)$
can be touched  on the initial local   interval
\([0,T_0]\).
For $j\geq1$, we set
\[
 r_j=2^j r_0,\qquad L_j=\log(e+r_j).
\]
Whenever the level $r_j$ is attained, let $t_j$ be its first
attainment time:
\[
 t_j=\inf\{t\geq0:\|\rho(t)\|_{\infty}=r_j\}.
\]
 We also put
\begin{equation*}
 I_j=[t_{j-1},t_j],\qquad
 E_j=\int_{I_j}D(t)\,\mathrm dt.
\end{equation*}
For each fixed terminal level $j$, we define the fixed
normalization
\be\la{ly5.1}
 H_j(t)^2=e+r_j^\eta+A(t)^2,\qquad
 Q_j=r_j^{-1/3}+\sum_{i=1}^j2^{i-j}E_i.
\ee
In particular, $H_j$ depends on the terminal index $j$, but its
definition is not changed as $t$ varies in $[0,t_j]$.
The intervals $I_i$ have disjoint interiors. The basic energy
inequality therefore gives
\[
 \sum_i E_i
 =\sum_i\int_{I_i}D(t)\,\mathrm dt
 \leq\int_0^S D(t)\,\mathrm dt\leq C.
\]
Consequently, changing the order of summation, we have
\begin{equation}\label{eq:summable}
 \begin{aligned}
 \sum_j Q_j
 &=\sum_j r_j^{-1/3}
   +\sum_j\sum_{i=1}^j2^{i-j}E_i\\
 &\leq \sum_j r_j^{-1/3}
   +\sum_i E_i\sum_{j=i}^{\infty}2^{i-j}
 \leq C.
 \end{aligned}
\end{equation}
We point out that all sums here and below are restricted to attained levels.

Recalling the definition
\[
 h \Div u=G+P,\qquad
 A^2=\int\left(\mu\omega^2+\frac{G^2}{h}\right)\dd,
 \qquad
 D=\int\left(\mu\omega^2+h(\Div u)^2\right)\dd,
\]
we use  the elementary inequalities
\[
 \frac{G^2}{h}\leq 2h(\Div u)^2+2\frac{P^2}{h},
 \qquad
 h(\Div u)^2\leq 2\frac{G^2}{h}+2\frac{P^2}{h},
\]
to   conclude
\[
 A^2\leq C(D+\|P\|_{ 2 }^2),\qquad
 D\leq C(A^2+\|P\|_{ 2 }^2).
\]
Furthermore,  the pressure estimate \eqref{eq:pressure2} tells that if $\rho\leq r_j$, then   $ \sup_{0\leq t\leq t_j}\|P(t)\|_{ 2 }^2
       \leq C_\varepsilon r_j^\varepsilon.$   Consequently,
\begin{equation}\label{eq:shift}
 \begin{gathered}
 D(t)\leq C H_j(t)^2\quad(0\leq t\leq t_j),\quad
 \int_0^{t_j}A(t)^2\,\mathrm dt
       \leq C_\varepsilon r_j^\varepsilon,\quad
 \int_0^{t_j}(D(t)+A(t)^2)\,\mathrm dt
       \leq C_\varepsilon r_j^\varepsilon,
 \end{gathered}
\end{equation}
where the constants   do not depend on the choice of
$r_0$.

The following lemma establishes a local logarithmic-energy estimate, which will be used in the next section to control the time integral of the effective viscous flux in \(L^\infty\).
\begin{lemma}\label{lem:log-energy}
For each attained level $r_j$, it holds
\begin{equation}\label{eq:local-ratio}
 F_j\triangleq\int_{I_j}\frac{B^2}{H_j^2}\dt
 \leq C_\tau r_j L_j Q_j.
\end{equation}
\end{lemma}

\begin{proof}
It follows from the momentum equation that
\begin{equation}\label{eqhhd987}
 \rho\dot u=\nabla G+\mu\nabla^{\perp}\omega.
\end{equation}
Therefore, on any time interval on which $\rho\leq r$, we get
\begin{equation}\label{eq:elliptic}
 \|\nabla G\|_{ 2 }+\|\nabla\omega\|_{ 2 }
       \leq C\sqrt r\,B,\qquad
 \|G\|_{ 2 }\leq C r^{\gamma/2}A.
\end{equation}
Multiplying \eqref{eqhhd987} by $2\dot{u}$ and integrating it over $\mathbb{R}^{2}$, we are able to show
\begin{align*}
 (A^2)'+2B^2
 &=\int\left(
 4G\det\nabla u
 -2\left(1+\frac{\rho }h\right)G(\Div u)^2
 -\mu\omega^2 \Div u\right)\,\mathrm dx\\
 &\quad+\int\left(
 \left(1+\frac{\rho }h\right)\frac{G^2}{h}\Div u
 +2\gamma\frac{PG}{h}\Div u\right)\dd.
\end{align*}
Noticing the following identity
\[
 G(\Div u)^2=\left(\frac{G^2}{h}+\frac{PG}{h}\right)\Div u,
\]
we find that the terms involving $G(\Div u)^2$ can be combined with the
last integral. Since
$\rho/h(\rho)\leq C$, the resulting coefficients
are uniformly bounded.

By the Hardy--BMO estimate (see \cite{coi1,coi3}), we have
\[
 \left|\int G\det\nabla u\,\mathrm dx\right|
 \leq C\|\nabla G\|_{ 2 }\|\nabla u\|_{ 2 }^2
 \leq C\sqrt r\,B D.
\]
The
Gagliardo--Nirenberg inequality gives
\[
 \left|\int\omega^2 \Div u\,\mathrm dx\right|
 \leq\|\omega\|_{ 4 }^2\|\Div u\|_{ 2 }
 \leq C\|\omega\|_{ 2 }\|\nabla\omega\|_{ 2 }\sqrt D
 \leq C\sqrt r\,B A\sqrt D.
\]
Due to $\|\Div u\|_{ 2 }\leq C\sqrt D$, the remaining cubic terms are bounded by
\[
 C\sqrt D\left(\left\|\frac{G^2}{h}\right\|_{ 2 }
                    +\left\|\frac{PG}{h}\right\|_{ 2 }\right).
\]
We thus conclude that
\begin{equation}\label{eq:first-energy}
 (A^2)'+2B^2
 \leq C\sqrt D\left(
       \left\|\frac{G^2}{h}\right\|_{ 2 }+
       \left\|\frac{PG}{h}\right\|_{ 2 }\right)
       +C\sqrt r\,B(A\sqrt D+D).
\end{equation}
For $0<\delta\leq1/4$, the sharp Gagliardo--Nirenberg inequality yields
\[
 \|G\|_{{2(1+\delta)/\delta}}^{1+\delta}
 \leq C\delta^{-(1+\delta)/2}
       \|G\|_{ 2 }^\delta\|\nabla G\|_{ 2 }.
\]
Notice the fact that
\[
 \delta^{-\delta/2}
 =\exp\left(\frac{\delta}{2}\log\frac1\delta\right)\leq C
 \quad(0<\delta\leq1/4),
\]
the factor $\delta^{-(1+\delta)/2}$ is bounded by
$C\delta^{-1/2}$. Also, by H\"older's inequality, one has
\[
 \begin{aligned}
 \left\|\frac{G^2}{h}\right\|_{ 2 }
 &=\left\||G|^{1+\delta}\frac{|G|^{1-\delta}}h\right\|_{ 2 }\\
 &\leq \|G\|_{{2(1+\delta)/\delta}}^{1+\delta}
       \left\|\frac{|G|^{1-\delta}}h
          \right\|_{{2/(1-\delta)}}\\
 &\leq C\|G\|_{{2(1+\delta)/\delta}}^{1+\delta}
       \left(\int\frac{G^2}{h}\dd
          \right)^{(1-\delta)/2}\leq C A^{1-\delta}
       \|G\|_{{2(1+\delta)/\delta}}^{1+\delta}.
 \end{aligned}
\]
Combining this estimate with \eqref{eq:elliptic} gives
\[
 \begin{aligned}
 \left\|\frac{G^2}{h}\right\|_{ 2 }
 \leq C\delta^{-1/2}
       A^{1-\delta}\|G\|_{ 2 }^\delta\|\nabla G\|_{ 2 } \leq C\delta^{-1/2}
       r^{(1+\gamma\delta)/2}AB.
 \end{aligned}
\]
We also use
\begin{equation}\label{eq:interpolation-quadratic}
 \left\|\frac{PG}{h}\right\|_{ 2 }
 \leq C\|P\|_{ 4 }
       \left\|\frac{G^2}{h}\right\|_{ 2 }^{1/2}.
\end{equation}
On $I_i$, we take $r=r_i$ and choose
\[
 \delta_i=\frac1{4L_i}.
\]
Then $0<\delta_i\leq1/4$ and
\[
 r_i^{\gamma\delta_i}
 =\exp\left(\frac{\gamma\log r_i}{4L_i}\right)\leq C_\gamma.
\]
It follows that
\begin{equation}\label{eq:quadratic-on-Ii}
 \left\|\frac{G^2}{h}\right\|_{ 2 }
       \leq C(r_iL_i)^{1/2}AB
       \qquad (t\in I_i).
\end{equation}
Inserting this in
\eqref{eq:interpolation-quadratic}, we get
\[
 C\sqrt D\left\|\frac{PG}{h}\right\|_{ 2 }
 \leq C(r_iL_i)^{1/4}\|P\|_{ 4 }
           D^{1/2}A^{1/2}B^{1/2}.
\]
By Young's inequality, we have for any
fixed sufficiently small $\kappa>0$ that
\[
 \begin{aligned}
 C(r_iL_i)^{1/4}\|P\|_{ 4 }
       D^{1/2}A^{1/2}B^{1/2}
  \leq \kappa B^2 +C_\kappa(r_iL_i)^{1/3}
          \|P\|_{ 4 }^{4/3}D^{2/3}(A^2)^{1/3},
 \end{aligned}
\]
which implies
\begin{equation}\label{eq:pressure-remainder-corrected}
 C\sqrt D\left\|\frac{PG}{h}\right\|_{ 2 }
 \leq \kappa B^2+
 C_\kappa(r_iL_i)^{1/3}\|P\|_{ 4 }^{4/3}(D+A^2).
\end{equation}
Throughout $[\tau,t_j]$, we use the same
normalization $H_j$ as defined in \eqref{ly5.1}. We thus get
\[
 (H_j^2)'=(A^2)',\qquad
 (\log H_j^2)'=\frac{(A^2)'}{H_j^2}.
\]
This along with \eqref{eq:first-energy} implies
\begin{equation}\label{eq:first-energy111}
 (\log H_j^2)'
       +2\frac{B^2}{H_j^2}
 \leq C\frac{\sqrt D}{H_j^2}\left(
       \left\|\frac{G^2}{h}\right\|_{ 2 }+
       \left\|\frac{PG}{h}\right\|_{ 2 }\right)
       +C\sqrt r_{i}\,\frac{B(A\sqrt D+D)}{H_j^2}.
\end{equation}

We estimate \eqref{eq:first-energy111} on each $I_i$ for $1\leq i\leq j$. By \eqref{eq:quadratic-on-Ii} and $\frac{A^2}{H_j^2}\leq1$, one has
\begin{align*}
 \frac{C\sqrt D\,\left\|\frac{G^2}{h}\right\|_2}{H_j^2}
 &\leq C(r_iL_i)^{1/2}
       \frac{AB\sqrt D}{H_j^2}
\\ &\leq \kappa\frac{B^2}{H_j^2}
       +C_\kappa r_iL_i
          \frac{A^2D}{H_j^2}
          \\ &\leq \kappa\frac{B^2}{H_j^2}
       +C_\kappa r_iL_iD.
\end{align*}
Dividing \eqref{eq:pressure-remainder-corrected} by
$H_j^2\geq1$, we get
\[
 \begin{aligned}
 \frac{C\sqrt D\,\left\|\frac{PG}{h}\right\|_{ 2 }}{H_j^2}
 &\leq \kappa\frac{B^2}{H_j^2} +C(r_iL_i)^{1/3}\|P\|_{ 4 }^{4/3}
                 \frac{D+A^2}{H_j^2}\\
 &\leq \kappa\frac{B^2}{H_j^2}
       +C(r_iL_i)^{1/3}\|P\|_{ 4 }^{4/3}(D+A^2).
 \end{aligned}
\]
The last term in \eqref{eq:first-energy111} can be bounded by
 \begin{align*}
 C\sqrt r_{i}\,\frac{B(A\sqrt D+D)}{H_j^2}&=C\sqrt{r_i}\frac{BA\sqrt D}{H_j^2}+
 C\sqrt{r_i}\frac{BD}{H_j^2} \\& \leq \kappa\frac{B^2}{H_j^2}
       +C_\kappa r_i\frac{A^2D}{H_j^2}+C_\kappa r_i\frac{D^2}{H_j^2}
       \\& \leq \kappa\frac{B^2}{H_j^2}+C_\kappa r_iD\\& \leq \kappa\frac{B^2}{H_j^2}+C_\kappa r_iL_{i}D,
\end{align*}
where we have used  $L_i\geq1$ and $D/H_j^2\leq C$ due to
\eqref{eq:shift}.
Choosing $\kappa$ so  small such that
\begin{equation}\label{eq:log-local}
 (\log H_j^2)'
       +c\frac{B^2}{H_j^2}
 \leq C r_iL_iD
       +C(r_iL_i)^{1/3}\|P\|_{ 4 }^{4/3}(D+A^2)
 \quad\hbox{on }I_i.
\end{equation}

We now integrate \eqref{eq:log-local} from $\tau$ to $t_j$.
Keeping in mind
$A(\tau)<\infty$ and $r_j\geq 1$, it follows that
\[
 \log H_j(\tau)^2
 =\log\bigl(e+r_j^\eta+A(\tau)^2\bigr)
 \leq \log\bigl(e+1+A(\tau)^2\bigr)+\eta\log r_j
 \leq C_\tau+\eta L_j.
\]
Thus the contribution of the
logarithmic derivative is bounded from below by
$-C_\tau-\eta L_j$. Here $C_\tau$ is finite and independent of $j$.
For the pressure contribution, note that $r_iL_i\leq r_jL_j$
for $i\leq j$. Hence, we derive
\[
 \begin{aligned}
 \sum_{i=1}^j
  (r_iL_i)^{1/3}
   \int_{I_i}\|P(t)\|_{ 4 }^{4/3}(D+A^2)\,\mathrm dt
 & \leq (r_jL_j)^{1/3}
   \left(\sup_{0\leq t\leq t_j}\|P(t)\|_{ 4 }^{4/3}\right)
   \int_0^{t_j}(D+A^2)\,\mathrm dt\\
 & \leq C_\varepsilon
       r_j^{1/3+2\varepsilon}L_j^{1/3},
 \end{aligned}
\]
where we used \eqref{eq:pressure2} and \eqref{eq:shift}. Because
$\varepsilon\leq1/12$ and $L_j^{1/3}\leq C r_j^{1/6}$, it is easy to check
\[
 r_j^{1/3+2\varepsilon}L_j^{1/3}
 \leq C r_j^{1/2}L_j^{1/3}
 \leq C r_j^{2/3}.
\]
It follows from \eqref{eq:log-local} that
\[
 \int_\tau^{t_j}\frac{B^2}{H_j^2}\,\mathrm dt
 \leq C_\tau\left(
      L_j+r_j^{2/3}
      +\sum_{i=1}^j r_iL_iE_i\right).
\]
Since $r_j\geq2$ and $L_j\geq1$, the term $L_j$ is absorbed by
$r_j^{2/3}$.
Because the integral over $I_j$ is no larger than the integral over
$[\tau,t_j]$, we are able to show
\[
 F_j\leq C_\tau\left(
      r_j^{2/3}+\sum_{i=1}^j r_iL_iE_i\right).
\]
Finally, the dyadic relation $r_i/r_j=2^{i-j}$ and monotonicity
$L_i\leq L_j$ imply
\[
 \sum_{i=1}^j r_iL_iE_i
 \leq r_jL_j\sum_{i=1}^j2^{i-j}E_i.
\]
Also, one has
\[
 r_j^{2/3}\leq r_jL_j\,r_j^{-1/3}.
\]
Recalling the definition of $Q_j$, we arrive at
\begin{equation*}
 F_j
 \leq C_\tau r_jL_jQ_j.
\end{equation*}
This ends the proof of Lemma \ref{lem:log-energy}.
\end{proof}

Finally, we use the time interval layer  foliation  to estimate the  effective  viscous flux  on each interval \(I_j\).

\begin{lemma}\label{lem:flux}
For each attained level, it holds
\begin{equation}\label{eq:flux}
 \int_{I_j}\frac{\|G\|_{ \infty }^2}{H_j^2}\dt
 \leq C_\tau r_j^2 L_j^2 Q_j.
\end{equation}
\end{lemma}

\begin{proof}
From the inequality (see Appendix for its proof)\begin{equation}\label{sdbnmygf}
(B^2)'+cW^2
 \leq  C\int
  \left((1+\lambda)^2|\nabla u|^4
        +P^2|\nabla u|^2\right)\dd,
\end{equation}
we see that,  if
$\rho\leq r_j$,
 \begin{align} \label{dfgD.1}
 (B^2)'+cW^2
 &\leq C r_j^{2\gamma}
       \left(\|\nabla u\|_4^4+\|\nabla u\|_2^2\right).
 \end{align}
The Gagliardo--Nirenberg
inequality gives
\begin{align} \label{dfgD.2}
 \|\nabla u\|_{ 4 }^4
 \leq C\bigl(\|\Div u\|_{ 4 }^4+\|\omega\|_{ 4 }^4\bigr)
 \leq C\bigl(\|\Div u\|_{ 4 }^4+
       \|\omega\|_{ 2 }^2\|\nabla\omega\|_{ 2 }^2\bigr).
 \end{align}
Since $h \Div u=G+P$, we decompose the domain into
\[
 \{|G|\geq \tfrac12 h|\Div u|\}
 \quad\hbox{and}\quad
 \{|G|<\tfrac12 h|\Div u|\}.
\]
On the first set, $|\Div u|\leq 2|G|/h$; on the second,
$|P|=|h \Div u-G|\geq h|\Div u|/2$, and hence
$|\Div u|\leq 2P/h$.  Using $h\geq2\mu$ and
$P/h\leq C r_j^{\gamma-1}$ on this set, this gives
\begin{align}
 |\Div u|^4\leq C|G|^4+C r_j^{2\gamma}|\Div u|^2.\nonumber
\end{align}
Consequently, we have
\begin{align} \label{dfgD.4}
 \|\Div u\|_4^4
 \leq C\|G\|_{ 4 }^4+C r_j^{2\gamma}\|\Div u\|_{ 2 }^2.
\end{align}
By the Gagliardo--Nirenberg, \eqref{eq:elliptic} and
$\|\Div u\|_2^2\leq C D$, we obtain
\begin{align} \label{dfgD.5}
 \|G\|_{ 4 }^4
 &\leq C\|G\|_{ 2 }^2\|\nabla G\|_{ 2 }^2
 \leq C r_j^\gamma A^2\,(r_jB^2)
 =C r_j^{\gamma+1}A^2B^2,\nonumber\\
 \|\omega\|_{ 4 }^4
 &\leq C\|\omega\|_{ 2 }^2\|\nabla\omega\|_{ 2 }^2
 \leq C r_j A^2B^2.
 \end{align}
Combining \eqref{dfgD.2}, \eqref{dfgD.4} and \eqref{dfgD.5}, we deduce by using $\|\Div u\|_2^2\leq C D$ that
 \begin{align}
 \|\nabla u\|_{ 4 }^4
 &\leq  C\bigl(\|\Div u\|_{ 4 }^4+
       \|\omega\|_{ 2 }^2\|\nabla\omega\|_{ 2 }^2\bigr)\nonumber\\
       &\leq  C\bigl(\|G\|_{ 4 }^4+C r_j^{2\gamma}\|\Div u\|_{ 2 }^2 +
       \|\omega\|_{ 2 }^2\|\nabla\omega\|_{ 2 }^2\bigr)\nonumber\\
       &\leq  C\bigl(r_j^{\gamma+1}A^2B^2+C r_j^{2\gamma}D+r_j A^2B^2\bigr)
       \nonumber\\
       &\leq  C\bigl(r_j^{\gamma+1}A^2B^2+C r_j^{2\gamma}D\bigr), \nonumber
 \end{align}
where we have used $r_j\geq2$.
We thus derive
\begin{align} \label{sapjkh28}
\|\nabla u\|_4^4+\|\nabla u\|_2^2\leq  C\bigl(r_j^{\gamma+1}A^2B^2+C r_j^{2\gamma}D\bigr).
 \end{align}
Putting \eqref{sapjkh28} into \eqref{dfgD.1} and using $\gamma>1$ as well as $r_j\geq2$, it follows that
\[
 \begin{aligned}
 (B^2)'+cW^2
 &\leq C\left(
 r_j^{3\gamma+1}A^2B^2+r_j^{4\gamma}D\right)\\
 &\leq C r_j^{4\gamma}(A^2B^2+D).
 \end{aligned}
\]
On $[0,t_j]$, one gets
\begin{equation}\label{eq:accel}
 (B^2)'+cW^2\leq C r_j^{4\gamma}(A^2B^2+D).
\end{equation}
Dividing \eqref{eq:accel} by $e+B^2$ and using the identities
\[
 \frac{(B^2)'}{e+B^2}
       =\bigl(\log(e+B^2)\bigr)',\qquad
 \frac{A^2B^2}{e+B^2}\leq A^2,\qquad
 \frac{D}{e+B^2}\leq D,
\]
we get by integration from $\tau$ to $t_j$ that
\[
 c\int_\tau^{t_j}\frac{W^2}{e+B^2}\,\mathrm{d}t
 \leq \log(e+B(\tau)^2)
       +C r_j^{4\gamma}\int_\tau^{t_j}(A^2+D)\,\mathrm{d}t.
\]
Since \(\tau\) is fixed and \(B(\tau)<\infty\), the initial logarithmic term is bounded by a constant depending only on \(\tau\).
By \eqref{eq:shift}, we arrive at
\begin{align} \label{dfgD.6}
 \int_\tau^{t_j}\frac{W^2}{e+B^2}\,\mathrm{d}t
 \leq C_\tau+C_\eps r_j^{4\gamma+\eps}
 \leq C_\tau r_j^{4\gamma+1}.
\end{align}
We have included one additional power of \(r_j\) in the last inequality to avoid tracking the harmless small power \(r_j^\varepsilon\).

Taking the divergence of the momentum equation gives
$
 \Div(\rho\dot u)=\Delta G.
$
Hence,
\[
 G=-I\Div(\rho\dot u)\quad{\rm with}\quad I=(-\Delta)^{-1}.
\]
We use the potential interpolation inequality
\cite[Lemma~3.2]{FJLL}:
\begin{align} \label{dfgD.7}
 \|G\|_{ \infty }
 \leq C\sqrt{r_j}\,B+
 C\sqrt{p r_j}\,
       \|G\|_{ 2 }^\sigma B^{1-2\sigma}W^\sigma,
\end{align}
where $\sigma=\frac{2}{p+2}$ with $p\geq4$.
We mention that the constant in this inequality is independent of $p$.
Moreover, the factor $\sqrt p$ is important: choosing $p$
proportional to $\log r_j$ will keep the density powers generated
by the interpolation uniformly bounded.

Recalling \eqref{eq:elliptic}, we have
\[
 \|G\|_2\leq C r_j^{\gamma/2}A.
\]
Squaring \eqref{dfgD.7}, using $(x+y)^2\leq2x^2+2y^2$ and dividing by
$H_j^2$, we get
\begin{align} \label{dfgD.8}
 \frac{\|G\|_{ \infty }^2}{H_j^2}
 \leq C r_j\frac{B^2}{H_j^2}
 +Cp r_j^{1+\gamma\sigma}
   \left(\frac{A^2}{H_j^2}\right)^\sigma
   \left(\frac{B^2}{H_j^2}\right)^{1-2\sigma}
   \left(\frac{W^2}{e+B^2}\right)^\sigma
   \left(\frac{e+B^2}{H_j^2}\right)^\sigma.
\end{align}
Since $H_j^2=e+r_j^\eta+A^2\geq e+1$, then
$
 \frac{e+B^2}{H_j^2}
 \leq 1+\frac{B^2}{H_j^2}.
$
Thus,  from \eqref{dfgD.8} we deduce
\begin{align} \label{dffv235}
 \frac{\|G\|_{ \infty }^2}{H_j^2}
 \leq &C r_j\frac{B^2}{H_j^2}
 +Cp r_j^{1+\gamma\sigma}
   \left(\frac{A^2}{H_j^2}\right)^\sigma
   \left(\frac{B^2}{H_j^2}\right)^{1-2\sigma}
   \left(\frac{W^2}{e+B^2}\right)^\sigma
  \nonumber\\&+
   Cp r_j^{1+\gamma\sigma}
   \left(\frac{A^2}{H_j^2}\right)^\sigma
   \left(\frac{B^2}{H_j^2}\right)^{1-\sigma}
   \left(\frac{W^2}{e+B^2}\right)^\sigma\nonumber\\
   \leq &C r_j\frac{B^2}{H_j^2}
 +Cp r_j^{1+\gamma\sigma}
   \left(\frac{A^2}{H_j^2}\right)^\sigma
   \left(\frac{B^2}{H_j^2}\right)^{1-2\sigma}
   \left(\frac{W^2}{e+B^2}\right)^\sigma
  \nonumber\\&+
   Cp r_j^{1+\gamma\sigma}
   \left(\frac{B^2}{H_j^2}\right)^{1-\sigma}
   \left(\frac{W^2}{e+B^2}\right)^\sigma,
\end{align}
where in the last line we have used $H_j^2\geq A^2$.
Denote by
\[
 X(t)=\frac{A(t)^2}{H_j(t)^2},\qquad
 Z(t)=\frac{W(t)^2}{e+B(t)^2}.
\]
It follows from \eqref{dffv235} that
\begin{align}
 \int_{I_j}\frac{\|G\|_{ \infty }^2}{H_j^2}\,\mathrm{d}t
 &\leq C r_jF_j +Cp r_j^{1+\gamma\sigma}
 \left[
  \left(\int_{I_j}X\,\mathrm{d}t\right)^\sigma
  F_j^{1-2\sigma}
  \left(\int_{I_j}Z\,\mathrm{d}t\right)^\sigma
  +F_j^{1-\sigma}
  \left(\int_{I_j}Z\,\mathrm{d}t\right)^\sigma
 \right],\nonumber
 \end{align}
where $F_j$ is given by \eqref{eq:local-ratio}.
By \eqref{eq:shift} and \eqref{dfgD.6}, one has
\[
 \int_{I_j}X\,\mathrm{d}t
 \leq\int_0^{t_j}A^2\,\mathrm{d}t
 \leq C_\eps r_j^\eps,
 \qquad
 \int_{I_j}Z\,\mathrm{d}t
 \leq C_\tau r_j^{4\gamma+1}.
\]
  Therefore, we achieve
\begin{equation}\label{eq:flux-holder}
 \begin{aligned}
 \int_{I_j}\frac{\|G\|_{ \infty }^2}{H_j^2}\,\mathrm{d}t
 &\leq C r_jF_j +C_\tau p\,
 r_j^{1+(5\gamma+1+\eps)\sigma}
       \bigl(F_j^{1-2\sigma}+F_j^{1-\sigma}\bigr).
 \end{aligned}
\end{equation}
We choose
\[
 p=8+L_j,\qquad \sigma=\frac{2}{p+2}.
\]
With this choice, we have $p\geq4$, $p\leq C L_j$, and
$\sigma\leq C/L_j$.  Since $L_j=\log(e+r_j)$, we derive
\begin{align} \label{dfgD.10}
 r_j^{(5\gamma+1+\eps)\sigma}
 =\exp\!\left((5\gamma+1+\eps)\sigma\log r_j\right)
 \leq C_{\gamma,\eps}.
\end{align}
According to \eqref{eq:local-ratio}, one has
\[
 F_j\leq C_\tau r_jL_jQ_j.
\]
Also, by definition $Q_j\geq r_j^{-1/3}$ and hence
$r_jL_jQ_j\geq r_j^{2/3}L_j\geq1$. Since
$0<1-2\sigma<1$ and $0<1-\sigma<1$, we get
\[
 F_j^{1-2\sigma}+F_j^{1-\sigma}
 \leq C_\tau r_jL_jQ_j.
\]
Combining this observation with \eqref{dfgD.10} and using $p\leq CL_j$, we conclude from
\eqref{eq:flux-holder} that
\begin{equation}
\displaystyle
 \int_{I_j}\frac{\|G(t)\|_{ \infty }^2}{H_j(t)^2}\,\mathrm{d}t
       \leq C_\tau r_j^2L_j^2Q_j.\nonumber
\end{equation}
This ends the proof of the lemma.
\end{proof}

\section{Uniform upper bound for the density}
In this section, we establish a uniform-in-time upper bound for the density, which is the key estimate needed to prevent mass concentration and to complete the global continuation argument.

\begin{proposition}\label{prop:density}
There holds
\be\la{okly5.1}
\sup_{0\leq t<\infty}\|\rho(t)\|_{ \infty }\leq C.\ee
\end{proposition}

\begin{proof}
Let \(X(t;s,x)\) be a particle trajectory.
Then,   the mass  equation gives
\[
 \frac{\mathrm d}{\mathrm dt}\rho(t,X(t;s,x))
       =-(\rho\,\Div u)(t,X(t;s,x)).
\]
Recalling the  definition of $\theta$ in \eqref{ly2.1}, one has
\begin{align} \label{dfgD.11}
 \frac{\mathrm d}{\mathrm dt}\theta(\rho(t,X(t;s,x)))
       =-(h(\rho)\Div u)(t,X(t;s,x))=-G(t,X(t;s,x))-P(t,X(t;s,x)).
\end{align}
We choose a point \(x_j\) such that \(\rho(t_j,x_j)=r_j\), and trace the corresponding particle trajectory backward in time. By continuity, there
is a last time $s_j<t_j$ at which this trajectory has density
$r_{j-1}=r_j/2$.  On   open set  $(s_j,t_j)$ its density
is strictly larger than $r_{j-1}$; otherwise it must cross $r_{j-1}$ again before it touches  $r_j$. This contradicts  the
choice of $s_{j}$, and thus,
\begin{align} \label{dfgD.12}
 \rho(t,X(t))\geq r_j/2
       \quad (s_j\leq t\leq t_j).
\end{align}
For $j=1$, the choice of $r_0$ larger than the density on
$[0,T_0]$ ensures $s_1\geq\tau$.  For $j\geq2$, the first
attainment of $r_j$ and the last-crossing construction imply
$s_j\geq t_{j-1}$. Hence, it holds
\[
 [s_j,t_j]\subset I_j.
\]
For simplicity, we denote $\Delta_j=t_j-s_j$.

Keeping in mind
\[
 \theta'(r)=\frac{h(r)}r
 =1+\frac{2\mu}{r},
\]
we have, for $r_j\geq2$,
\[
 \begin{aligned}
 \theta(r_j)-\theta(r_{j-1})
 =\int_{r_j/2}^{r_j}
       \left(1+\frac{2\mu}{r} \right)
       \,\mathrm{d}r \geq \int_{r_j/2}^{r_j} \,\mathrm{d}r \geq c r_j.
 \end{aligned}
\]
Integrating \eqref{dfgD.11} along the selected trajectory gives
\begin{align} \label{sdgbkp1}
 \theta(r_j)-\theta(r_{j-1})
 =-\int_{s_j}^{t_j}G(t,X(t))\,\mathrm{d}t
  -\int_{s_j}^{t_j}P(t,X(t))\,\mathrm{d}t.
\end{align}
Taking the absolute value only in the first term at the right hand-side of \eqref{sdgbkp1}, and using
\eqref{dfgD.12}, we obtain
\begin{align} \label{dfgD.14}
 c r_j
 \leq\int_{s_j}^{t_j}\|G(t)\|_{ \infty }\,\mathrm{d}t
       -2^{-\gamma}r_j^\gamma\Delta_j.
\end{align}
We remark  that we do not take the absolute value of the pressure, since its sign is essential. In fact,
$P=\rho^\gamma\geq(r_j/2)^\gamma$ along this trajectory.

Using the
Cauchy--Schwarz inequality, we have
\begin{align} \label{dfgD.15}
 \int_{s_j}^{t_j}\|G(t)\|_{ \infty }\,\mathrm{d}t
 &\leq
 \left(\int_{s_j}^{t_j}
       \frac{\|G(t)\|_{ \infty }^2}{H_j(t)^2}\,\mathrm{d}t\right)^{1/2}
 \left(\int_{s_j}^{t_j}H_j(t)^2\,\mathrm{d}t\right)^{1/2}\nonumber\\
 &\leq C_\tau r_jL_j\sqrt{Q_j}
       \left(\int_{s_j}^{t_j}H_j(t)^2\,\mathrm{d}t\right)^{1/2},
 \end{align}
where we used \eqref{eq:flux} and
$[s_j,t_j]\subset I_j$.
 Since
$H_j^2=e+r_j^\eta+A^2$, it follows from \eqref{eq:shift} that
\[
 A^2\leq C(D+\|P\|_2^2),\qquad
 \|P(t)\|_2^2\leq C_\eps r_j^\eps
       \leq C_\eps r_j^\eta
       \quad(0\leq t\leq t_j),
\]
which implies
\[
 H_j(t)^2\leq C_\eps\bigl(D(t)+r_j^\eta\bigr).
\]
Consequently, we get
\begin{align} \label{dfgD.16}
 \int_{s_j}^{t_j}H_j(t)^2\,\mathrm{d}t
 \leq C_\eps\left(
       \int_{I_j}D(t)\,\mathrm{d}t+r_j^\eta\Delta_j\right)
 =C_\eps(E_j+r_j^\eta\Delta_j).
\end{align}
Combining \eqref{dfgD.14}, \eqref{dfgD.15} with  \eqref{dfgD.16}, we deduce
\be\ba\label{eq:cross}
 c r_j
 &\leq C_\tau r_jL_j\sqrt{Q_j}
       \sqrt{E_j+r_j^\eta\Delta_j}
       -2^{-\gamma}r_j^\gamma\Delta_j\\
       &\leq C_\tau r_jL_j\sqrt{Q_jE_j} +C_\tau r_j^{1+\eta/2}L_j\sqrt{Q_j}\sqrt{\Delta_j}
       -2^{-\gamma}r_j^\gamma\Delta_j,
\ea\ee
where the second inequality follows from the fact $\sqrt{a+b}\leq\sqrt a+\sqrt b$.    Notice that
\[
 az-bz^2\leq\frac{a^2}{4b},\quad \forall a,b>0,z\geq0,
\]
we can  choose in \eqref{eq:cross}
$a=C_\tau r_j^{1+\eta/2}L_j\sqrt{Q_j}$,
$b=2^{-\gamma}r_j^\gamma$ and $z=\sqrt{\Delta_j}$, to  conclude that
\begin{equation}\label{eq:cross-max}
 c r_j
 \leq C_\tau r_jL_j\sqrt{E_jQ_j}
       +C_\tau r_j^{2+\eta-\gamma}L_j^2Q_j.
 \end{equation}

 We next prove
 \be\la{ly5.2}
 C_\tau r_j^{\,1+\eta-\gamma}L_j^{2}
       \leq c/2\quad{\rm for\,\,every}\,\,j=1,2,....
\ee
to  guarantee that  the  last  term in \eqref{eq:cross-max}  can be absorbed  by the  left
side.  In fact,   \eqref{eq:summable} implies $Q_j\leq C$, and
therefore,  the last term    $C_\tau r_j^{2+\eta-\gamma}L_j^2Q_j \le C_\tau r_j^{\,1+\eta-\gamma}L_j^{2}\cdot  r_{j}$.
Since  the exponent $1+\eta-\gamma$ is strictly
negative, we see that   $r_j^{\,1+\eta-\gamma}=2^{j(1+\eta-\gamma)} r_0^{1+\eta-\gamma}$ is exponential decay in $j$ and   dominates
the polynomial growth  $L_j\asymp j+\log r_0$.
Therefore, \eqref{ly5.2} holds true, as long as  we choose  $r_0$  large enough (without changing the
previous estimates).

In terms of \eqref{eq:cross-max} and \eqref{ly5.2},  we obtain
 \begin{equation}\label{eq:level-cost}
  L_j^{-1}\leq C_\tau\sqrt{E_jQ_j}.
  \end{equation}
Suppose that the levels $r_1,\ldots,r_N$ are all attained.
Summing \eqref{eq:level-cost} and applying Cauchy--Schwarz gives
 \begin{align}\label{dfgD.18}
 \sum_{j=1}^N L_j^{-1}
  \leq C_\tau\sum_{j=1}^N\sqrt{E_jQ_j} \leq C_\tau
       \left(\sum_{j=1}^N E_j\right)^{1/2}
       \left(\sum_{j=1}^N Q_j\right)^{1/2} \leq C_\tau.
 \end{align}
The last line follows from the disjointness of the intervals
$I_j$, the basic dissipation estimate
$\sum_jE_j\leq\int_0^S D\,\mathrm{d}t\leq C$, and
the summability \eqref{eq:summable} of $Q_j$.

On the other hand, $r_j=2^jr_0$ implies
$L_j=\log(e+2^jr_0)\asymp j+1$. Hence, we obtain
\[
 \sum_{j=1}^\infty L_j^{-1}
 \asymp \sum_{j=1}^\infty (j+1)^{-1}
 =\infty.
\]
This contradicts the uniform
bound \eqref{dfgD.18} if arbitrarily many levels are attained.

It follows that the number of attained dyadic levels is bounded
by a constant independent of the compact lifespan $[0,S]$.
Since $r_j=2^jr_0$, the density cannot exceed the last attained
level.  We have therefore proved
\[
 \sup_{0\leq t\leq S}\|\rho(t)\|_{\infty(\R^2)}
 \leq C,
\]
where $C$ is independent of $S$. \end{proof}

\section{Proof of Theorem \ref{TH11}}\label{sec:continuation}

Once the upper bound for the density \eqref{okly5.1} is obtained, an argument similar to that in
\cite{hlia,huang-li-cauchy-published} shows that,  for any $T<\infty$,
\bnn
 \sup_{0<t\le T}
 \left[
 \|\bar x^a\rho(\cdot,t)\|_{L^1\cap H^1\cap W^{1,q}}
 +\|\nabla u(\cdot,t)\|_{2}
 +\|\sqrt{\rho}u(\cdot,t)\|_{2}
 \right]\leq C_{T}
\enn
and that, for some large number $N_{*}$ which depends on $T$,
\bnn
 \int_{B_{N_{*}}}\rho(T,x)\dd
 \geq \frac{1}{8} \int \rho_0(x)\dd.\enn
At this stage, we may restart the local existence argument at a time $T_{0}<T$ chosen sufficiently close to $T$. Then, by  Lemma \ref{lem2.1} we can extend the unique strong  solution beyond the  $(0,T]$.
The proof of Theorem \ref{TH11} is thus completed.

\appendix
\section{Proof of  \eqref{sdbnmygf}.}\label{sec:energy-tests}

Setting
\[
\begin{aligned}
 S_{ij} \triangleq\mu\partial_j u_i+\bigl((\mu+\lambda)\Div u-P\bigr)\delta_{ij},\quad
 \Sigma_{ij} \triangleq D_t S_{ij}+\Div u S_{ij}-(\partial_\ell u_j)S_{i\ell}.
\end{aligned}
\]
We  claim  that
\begin{equation}\label{flux:Qasd231}
 \left|\Sigma-\mu\nabla\dot u
       -(\mu+\lambda)(\Div \dot{u})I\right|
 \leq C\left((1+\lambda)|\nabla u|^2
             +P|\nabla u|\right).
\end{equation}
In fact, direct computation yields
\begin{equation}  \label{yesaq92}
  D_t(\partial_ju_i)
  =\partial_j\dot u_i
   -(\partial_ju_\ell)\partial_\ell u_i.
\end{equation}
Taking the trace (with the appropriate indices) yields
\begin{equation}\label{yesaq93}
  D_t(\Div u)
  =D_t(\partial_ju_j)=\partial_j\dot u_j
    -(\partial_ju_\ell)\partial_\ell u_j=\Div \dot{u}
    -(\partial_ju_\ell)\partial_\ell u_j.
\end{equation}
It follows from mass  equation that
\begin{equation}
  D_t\lambda=-\lambda \Div u,
  \qquad
  D_tP=-\gamma P \Div u.\nonumber
\end{equation}
This and  \eqref{yesaq93}  provide us
\[
\begin{aligned}
D_t\big((\mu+\lambda)\Div u-P\big)
&=(\mu+\lambda)D_t\Div u+(D_t\lambda)\Div u-D_tP\\
&=(\mu+\lambda)\Div \dot{u}-(\mu+\lambda)
       (\partial_ju_\ell)\partial_\ell u_j
       - \lambda (\Div u)^2+\gamma P \Div u.
\end{aligned}
\]
Making use of the identity
\[
  D_tS_{ij}
  =\mu D_t(\partial_ju_i)
   +D_t\big((\mu+\lambda)\Div u-P\big)\delta_{ij}
\]
and \eqref{yesaq92}, we obtain
\begin{equation}\label{yesaq94}
\begin{aligned}
D_tS_{ij}
= &\mu\partial_j\dot u_i
  +(\mu+\lambda)(\Div \dot{u})\delta_{ij} -\mu(\partial_ju_\ell)\partial_\ell u_i\\
&-\Big((\mu+\lambda)
       (\partial_ju_\ell)\partial_\ell u_j
       + \lambda (\Div u)^2-\gamma P \Div u\Big)\delta_{ij}.
\end{aligned}
\end{equation}
Moreover, we have
\begin{equation}\label{yesaq95}
\begin{aligned}
\Div uS_{ij}
&=\mu \Div u\,\partial_ju_i
  +\big((\mu+\lambda)(\Div u)^2-P\Div u\big)\delta_{ij},
\end{aligned}
\end{equation}
\begin{equation}\label{yesaq96}
\begin{aligned}
-(\partial_\ell u_j)S_{i\ell}
= -\mu(\partial_\ell u_j)(\partial_\ell u_i) -\big((\mu+\lambda)\Div u-P\big)\partial_i u_j.
\end{aligned}
\end{equation}
Combining \eqref{yesaq94}, \eqref{yesaq95} with
\eqref{yesaq96}, we find
\[
  \Sigma_{ij}
  =\mu\partial_j\dot u_i
   +(\mu+\lambda)(\Div \dot{u})\delta_{ij}
   +\widetilde{R}_{ij},
\]
where the remainder is given by
\begin{equation}\label{yesaq97}
\begin{aligned}
\widetilde{R}_{ij}
= &-\mu(\partial_ju_\ell)\partial_\ell u_i
    -\mu(\partial_\ell u_j)\partial_\ell u_i-\big((\mu+\lambda)\Div u-P\big)\partial_i u_j
    +\mu \Div u\,\partial_ju_i\\
&+\Big[
   \mu (\Div u)^2
   -(\mu+\lambda)(\partial_ju_\ell)\partial_\ell u_j
   +(\gamma-1)P\Div u
  \Big]\delta_{ij}.
\end{aligned}
\end{equation}
In particular, it holds
\[
  \widetilde{R}=\Sigma-\mu\nabla\dot u
       -(\mu+\lambda)(\Div \dot{u})I.
\]
The terms in \eqref{yesaq97} can be bounded by
\[
 \left|
-\mu(\partial_ju_\ell)\partial_\ell u_i
-\mu(\partial_\ell u_j)\partial_\ell u_i
+\mu \Div u\,\partial_ju_i
\right| \leq C\mu|\nabla u|^2,
\]
\[
 \left|
\mu (\Div u)^2
-(\mu+\lambda)(\partial_ju_\ell)\partial_\ell u_j
\right|\leq C(1+\lambda)|\nabla u|^2,
\]
\[
  |(\gamma-1)P\Div u|
  \leq CP|\nabla u|,
\]
\[
\left|
\big((\mu+\lambda)\Div u-P\big)\partial_i u_j
\right|\leq
  (\mu+\lambda)|\Div u|\,|\nabla u|
  +P|\nabla u|\leq
  C(1+\lambda)|\nabla u|^2
  +P|\nabla u|.
\]
Combining these estimates, we obtain
\[
  |\widetilde{R}|
  \leq C\left((1+\lambda)|\nabla u|^2
             +P|\nabla u|\right).
\]
Recalling the definition of $\widetilde{R}$, we thus conclude  \eqref{flux:Qasd231}.

We next show that
\begin{equation}  \label{yesaq98}
 \rho D_t\dot u_i=\partial_j\Sigma_{ij}.
\end{equation}
To see this,  a straightforward computation shows that
\begin{equation}  \label{ssa235}
  D_t(\partial_jS_{ij})
  =\partial_j(D_tS_{ij})
   -(\partial_j u_\ell)\partial_\ell S_{ij}.
\end{equation}
and
\[
\begin{aligned}
\partial_j\Sigma_{ij}
={}&\partial_j(D_tS_{ij})
  +\partial_j(\Div u S_{ij})
  -\partial_j\big((\partial_\ell u_j)S_{i\ell}\big)\\
={}&\partial_j(D_tS_{ij})
  +(\partial_j\Div u)S_{ij}
  +\Div u\,\partial_jS_{ij}-(\partial_j\partial_\ell u_j)S_{i\ell}
  -(\partial_\ell u_j)\partial_jS_{i\ell}.
\end{aligned}
\]
Since $\partial_j\partial_\ell u_j=\partial_\ell \Div u$,
the two terms containing a derivative of $\Div u$ cancel after relabeling
the dummy index. Consequently,
\[
\begin{aligned}
\partial_j\Sigma_{ij}
 =\partial_j(D_tS_{ij})
  +\Div u\,\partial_jS_{ij}
  -(\partial_\ell u_j)\partial_jS_{i\ell} =D_t(\partial_jS_{ij})+\Div u\,\partial_jS_{ij},
\end{aligned}
\]
where the last equality follows from \eqref{ssa235}. Thus we derive
\begin{equation*}
  \partial_j\Sigma_{ij}
  =D_t(\partial_jS_{ij})+\Div u\,\partial_jS_{ij}.
\end{equation*}
By $\eqref{eq:system}_{2}$, we get
\[
  \partial_jS_{ij}=\rho\dot u_i,
\]
which together  with the mass  equation gives  that
\[
\begin{aligned}
\partial_j\Sigma_{ij}
  =D_t(\rho\dot u_i)+\Div u\rho\dot u_i =(D_t\rho)\dot u_i+\rho D_t\dot u_i
   +\Div u\rho\dot u_i =\rho D_t\dot u_i.
\end{aligned}
\]
Therefore,   \eqref{yesaq98} is valid.

Multiplying \eqref{yesaq98} by $\dot u_i$ and integrating
over $\mathbb R^2$, we arrive at
$$ \int \rho\dot u_iD_t\dot u_i\dd= \int \dot u_i\partial_j\Sigma_{ij}\dd.$$
By mass  equation, the left-hand side satisfies
\[
\begin{aligned}
 \int \rho\dot u_iD_t\dot u_i\dd
&= \int
  \rho\left(\partial_t\dot u_i
       +u_j\partial_j\dot u_i\right)\dot u_i\dd\\
&=\frac12 \int \rho\,
       \partial_t|\dot u|^2\dd
  +\frac12 \int \rho u\cdot
       \nabla|\dot u|^2\dd\\
&=\frac12 \int \rho\,
       \partial_t|\dot u|^2\dd
  -\frac12 \int
       \operatorname{div}(\rho u)|\dot u|^2\dd
=\frac12\frac{d}{dt}
        \int \rho|\dot u|^2\dd.
\end{aligned}
\]
For the right-hand side, integration by parts gives
\[
 \int \dot u_i\partial_j\Sigma_{ij}\dd
 =- \int
       \Sigma_{ij}\partial_j\dot u_i\dd =- \int
       \Sigma:\nabla\dot u\dd.
\]
Combining the two identities, and recalling the definition of $B^2$,
we conclude that
\[
  \frac12\frac{d}{dt}B^2(t)=  \frac12\frac{d}{dt}
   \int \rho|\dot u|^2\dd
  =- \int \Sigma:\nabla\dot u\dd.
\]
By Young's inequality,
\[
\begin{aligned}
\left| \int \widetilde{R}:\nabla\dot u\,dx\right|
\leq \frac{\mu}{2}\|\nabla\dot u\|_{2}^2 +C \int
  \left((1+\lambda)^2|\nabla u|^4
        +P^2|\nabla u|^2\right)\dd.
\end{aligned}
\]
Thus, using \eqref{flux:Qasd231}, we obtain
\[
\begin{aligned}
\frac{d}{dt}B^2(t)
&+\mu\|\nabla\dot u\|_{2}^2
 +2 \int (\mu+\lambda)
       (\Div \dot{u})^2\dd\\
&\leq C \int
  \left((1+\lambda)^2|\nabla u|^4
        +P^2|\nabla u|^2\right)\dd,
\end{aligned}
\]
which along with the definition of $W=\|\nabla\dot u\|_{2}$ yield the desired  \eqref{sdbnmygf}.

\vskip .4in
\textbf{Statements and Declarations:} The authors declare that there is no conflict of interest.
\vskip .1in

\textbf{Data availability:} Since this work is of abstract theoretical nature, no data sets are generated or analyzed.

\vskip .1in
\textbf{Declaration of AI Use:}
The authors acknowledge the use of AI tools for exploratory discussions and language refinement during the preparation of this work. All mathematical statements, proofs, and references were independently checked, revised, and finalized by the authors, who take full responsibility for the content of this paper.

\vskip .1in \noindent \textbf{Acknowledgments}\,\,
Liang was partially supported by the National Natural Science Foundation of China (Grant No. 12531010) and the Sichuan Science and Technology Program (Grant No. 2025ZYD0157). Ye was supported by the National Natural Science Foundation of China (Grant No. 12671274).
The  authors would like to thank Professor Jing Li for his interest in this work and for his valuable suggestions.

\end{document}